\documentclass[final]{amsart}
\author{D. Peter Overholser}
\address{KU Leuven\\
Celestijnenlaan 200B\\
3001 Leuven, Belgium}
\title{Descendent Tropical Mirror Symmetry for $\mathbb{P}^2$}

\usepackage{amssymb}
\usepackage{tikz}
\usetikzlibrary{arrows}
\usetikzlibrary{decorations.markings}
\tikzstyle arrowstyle=[scale=1]
\usepackage{pgfplots}
\pgfplotsset{compat=newest}
\usetikzlibrary{calc}
\tikzstyle dir=[postaction={decorate,decoration={markings,
    mark=at position .65 with {\arrow[arrowstyle]{stealth}}}}]
\tikzstyle rdir=[postaction={decorate,decoration={markings,
    mark=at position .65 with {\arrowreversed[arrowstyle]{stealth};}}}]
\usepackage{graphicx}
\usepackage{amsthm}
\usepackage{amsrefs}
\usepackage{amsmath}
\usepackage{subfig}
\usepackage{graphicx}
\newcommand {\ZZ} {\mathbb{Z}}
\newcommand {\RR} {\mathbb{R}}

\newcommand {\Spec} {\operatorname{Spec}}

\newcommand {\Hom} {\operatorname{Hom}}
\newcommand\diff{\operatorname{diff}}

\newcommand\gr{\operatorname{gr}}
\newcommand\op{\operatorname{op}}
\newcommand\tilop{\operatorname{\tilde{op}}}

\newtheorem{theorem}{Theorem}[section]

\newtheorem{thm}{Theorem}[section]
\newtheorem{lemma}[theorem]{Lemma}
\newtheorem{sublemma}[theorem]{Sublemma}
\newtheorem{claim}[theorem]{Claim}
\newtheorem{proposition}[theorem]{Proposition}
\newtheorem{corollary}[theorem]{Corollary}

\theoremstyle{definition}
\newtheorem{definition}[theorem]{Definition}
\newtheorem{definition-lemma}[theorem]{Definition-Lemma}

\theoremstyle{remark}

\numberwithin{equation}{section}
\numberwithin{figure}{section}

\begin{document}
\begin{abstract}
We modify Gross's construction of mirror symmetry for $\mathbb{P}^2$ \cite{MSP2} by introducing a \emph{descendent} tropical Landau-Ginzburg potential.  The period integrals of this potential compute a modification of Givental's $J$-function, explicitly encoding a larger sector of the big phase space.  As a byproduct of this construction, new tropical methods for computing certain descendent Gromov-Witten invariants are defined.
\end{abstract}
\maketitle

\section{Introduction}
Following the pioneering work of Mikhalkin \cite{Mik}, tropical geometry has become a significant force in enumerative geometry and mirror symmetry.   The Strominger-Yau-Zaslow conjecture \cite{SYZ} provides a helpful heuristic for understanding this success, positing that mirror manifolds $X$ and $\check{X}$ possess dual special Lagrangian fibrations over a common base $B$ (see \cite{gsyz} for a discussion of some recent advances in this direction).  The natural geometry on $B$ is tropical, and certain predictions of mirror symmetry are expected to be apparent as identifications of tropical structures that regulate both the $A$-model of $X$ and the $B$-model of $\check{X}$.  

Such an identification has been made in Gross's study of mirror symmetry for $X:=\mathbb{P}^2$ \cite{kan} \cite{MSP2}, which serves as a framework for this paper.  As introduced by Givental in \cite{homgiv}, we will consider a relationship between $X$ and a mirror \emph{Landau-Ginzburg (LG) model}.  An LG model is a pair $(\check{X}, W)$, where $\check{X}$ is a manifold (in this case $(\mathbb{C}^*)^2$) and $W:\check{X}\rightarrow \mathbb{C}$ is a regular function.  This work was expanded upon by Barannikov \cite{bar}, who used period integrals and semi-infinite variation of Hodge parameters to construct a Frobenius manifold structure (see \cite{man}) on the (B-model) moduli, the \emph{univeral unfolding} of the LG potential.  Under the mirror map, this Frobenius manifold is identified with another, that defined by the big quantum cohomology ring (A-model) of $X$.  The induced change of coordinates relates the moduli parameters of the LG model with the flat coordinates of quantum cohomology.  In these coordinates, the period integrals of the $B$-model compute Givental's $J$-function.

Work of Mikhalkin \cite{Mik} established the equality of certain tropical curve counts in $\mathbb{R}^2$ with Gromov-Witten (henceforth GW) invariants of $\mathbb{P}^2$, giving a tropical interpretation of its A-model.  A corresponding tropical structure on the mirror side was introduced by Gross.  One motivational principle was Cho and Oh's \cite{CO} work relating terms in the LG potential to counts of Maslov index 2 holomorphic disks in $X$.  See also \cite{aur}. A related construction was given by Fukaya, Oh, Ohta, and Ono in \cite{FOOO}, but in a category in which concrete calculation is quite difficult.  By reducing to the tropical setting, one can define a combinatorial analogue of holomorphic disks (as shown in the work of Nishinou \cite{N}), which in turn can be used to define a \emph{tropical} unfolding of the LG potential.  Under the period integral, these tropical disks are glued together into curves whose related GW invariants govern the quantum cohomology, and the mirror symmetric coordinate transformation is extremely natural.  It is important to note, however, that mirror symmetry for $\mathbb{P}^2$ (in contrast to the elliptic curve case \cite{BM}) relies on GW invariants that have no \emph{a priori} tropical interpretation.  The mirror map is used to give such an interpretation as counts of objects assembled from tropical disks; the foregoing argument does not yield a tropical proof of this type of mirror symmetry for $\mathbb{P}^2$, but instead demonstrates an equivalence between the classical validity of certain tropical GW invariants and mirror symmetry.  

We will pursue a variation on this theme.  Instead of studying tropical structures that are expected to relate to classical objects of interest, we will work with combinatorial constructions of unknown classical significance and study their behavior under mirror symmetry.  We define a set of \emph{descendent} tropical disks that result in a finer unfolding of the LG potential.  The period integrals lead to a tropical version $\mathbb{T}_{trop}$ of a modification of Givental's $J$-function that explicitly encodes a larger than usual sector of the big phase space.  Mirror symmetry produces a very nicely behaved change of coordinates, under which $\mathbb{T}_{trop}$ is identified with a pullback $\mathbb{J}$ of the $J$-function.  We then use axioms of GW theory to show that $\mathbb{T}_{trop}$ is equal to its classical counterpart $\mathbb{T}$, a descendent version of the $J$-function.   

\subsection{Results}
We paraphrase our results here.  \begin{enumerate}
\item  A set of points in general position in the plane determines a \emph{tropical descendent unfolding} of the LG potential $W$.  There exists additional data $\Xi_i$, $f$, and $\Omega$ such that the corresponding period integrals satisfy the conditions of Section 1 of \cite{MSP2} and do not depend on the choice of general points.   

\item  The period integrals calculate a generating function $\mathbb{T}_{trop}$ of tropical invariants of type 
$$\langle \psi^\nu T_{i,trop}, T_{0,trop}^k, \psi^{r_1} T_{2,trop}, \ldots, \psi^{r_n} T_{2,trop}\rangle^{trop}_{0,d},$$
(degree $d\geq 0$, genus $0$) where $T_{i,trop}$ is a tropical analogue of the condition given by a positive generator $T_i$ of $H^{2i}(\mathbb{P}^2, \mathbb{Z})$.  We will later justify the suggestive notation. \label{item1}

\item The flat coordinates on the B-model Frobenius manifold are given as a generating function of tropical curve of type
$$\langle  T_{i,trop}, T_{0,trop}^k, \psi^{r_1} T_{2,trop}, \ldots, \psi^{r_n} T_{2,trop}\rangle^{trop}_{0,d}.$$
By the results of \cite{M+R}, these invariants are equal to their classical counterparts, and thus the change of coordinates induced by mirror symmetry is defined in terms of classical GW theory.\label{item2}

\item We define $\mathbb{J}$ as the pullback of the $J$-function to the B-model Frobenius manifold and $\mathbb{T}$ to be the classical analogue of $\mathbb{T}_{trop}$.  By \cite{bar}, $\mathbb{J}=\mathbb{T}_{trop}$.\label{item3}

\item $\mathbb{T}$ can be condensed as a generating function of GW invariants of type
$$\langle \psi^\nu T_{i}, T_{0}^k, T_{1}^l,\psi^{r_1} T_{2}, \ldots, \psi^{r_n} T_{2}\rangle_{0,d}.$$ \label{item4}

\item $\mathbb{J}=\mathbb{T}$, and thus the tropical curve counts of Item \ref{item1} are equal to their corresponding classical GW invariants.  Therefore, the period integrals compute a generating function encoding all classical GW invariants of the type given in Item \ref{item4}. \label{item5}

\end{enumerate}

\subsection{Overview}
Section \ref{def} is a set of preliminary definitions.  The reader is advised to pay close attention to $\mathcal{R}_k$, a regrettably intricate bookkeeping structure with some unusual operations.  

In Section \ref{trop}, we define the tropical objects necessary for our construction, which are purely combinatorial in nature.  Our tropical curves are generalizations of those found in \cite{kan} and \cite{M+R}, for they are designed to calculate genus zero invariants of type  $$\langle \psi^\nu T_i, T_0^k, \psi^{r_1} T_2, \ldots, \psi^{r_n} T_2\rangle_{0,d}.$$  In \cite{kan}, Gross gives tropical methods to calculate invariants of type $$\langle \psi^\nu T_i,  T_2, \ldots,T_2\rangle_{0,d},$$ while Markwig and Rau \cite{M+R} use an intersection theory to give tropical versions of $$\langle T_0^k, T_1^l, \psi^{r_1} T_2, \ldots, \psi^{r_n} T_2\rangle_{0,d}.$$  From a combinatorial perspective, the addition of a $\psi$ class in a descendent tropical invariant is reflected by an increment in the required valency of a vertex in its corresponding tropical curves.
 
We will also make use of a modification of the concept of the \emph{tropical disk} found in \cite{MSP2}; these should be understood as fragments of tropical curves broken apart at a vertex.  Instead of restricting to trivalent disks, we will allow higher valence vertices to occur at marked points.  The valences are recorded using $\mathcal{R}_k$

We assemble our tropical objects into moduli spaces of predictable dimension and define certain counts of tropical curves as putative tropical descendent GW invariants.  Their invariance and relation to classical GW theory is justified in later sections.

In Section \ref{bmod}, we introduce the B-model moduli relevant to our problem.  In \cite{kan}, the tropical LG potential is defined as a sum of monomials defined by trivalent disks passing through a selection of points in the plane, while the sum for our descendent potential runs over disks with higher valence vertices.  

The period integrals of \cite{kan} are adapted to this setting in Section \ref{integrals}, which is condensed from \cite{thesis}.  The process involves a generalization of the \emph{scattering diagrams} and \emph{broken lines} of \cite{MSP2}, accommodating the presence of higher valence vertices.  
We show that our descendent LG potential satisfies certain wall-crossing rules, and use this to prove that the integrals do not depend on the choice of general points used to determine the potential.  Furthermore, the period integrals extract a generating function whose coefficients are the descendent tropical invariants referenced above in Item \ref{item1}.  See Theorem \ref{theorem1}.

Section \ref{form} treats a number of formal manipulations of generating series of tropical and classical descendent GW invariants. We ``normalize" the integrals to satisfy the requirements of Section 1 of  \cite{MSP2}, allowing us to apply mirror symmetry.  The generating function $\mathbb{T}_{trop}$ defined by the integrals can then be related to a pullback $\mathbb{J}$ of the $J$-function by identifying flat coordinates in the B-model moduli.   See Theorem \ref{theorem2}.  These flat coordinates can be written in terms of classical GW invariants, yielding an expression for $\mathbb{T}_{trop}$ entirely in these terms.  Finally, the axioms of GW theory are applied to show that $\mathbb{J}$ is equal to the classical counterpart $\mathbb{T}$ of $\mathbb{T}_{trop}$.  This result confirms the classical relevance of the tropical descendent LG potential and the descendent GW invariants defined in Section \ref{trop}.
\subsection{Remarks}
There are several natural directions for further study.  Most immediately, it seems clear that a careful variation on the techniques above would compute tropical (and classical) GW invariants with arbitrary divisor insertions.  It is less clear how one could modify the techniques to allow further descendent insertions of non-$T_2$ classes.  In another vein, the techniques of this paper could be generalized to $\mathbb{P}^n$ for $n>2$.

 As in all approaches to tropical curve counting, the concept of multiplicity plays a central role here.  Mikhalkin's famous multiplicity \cite{Mik}, central to this paper, is now well understood even in higher dimensions \cite{NS}.  In contrast, some of the other multiplicities encountered here and in \cite{MSP2} are still mysterious, but may help to build stronger connections between the tropical and classical world.  In particular, they should be seen as manifestations of a potential log-geometric construction for $\mathbb{P}^2$ analogous to that worked out for $\mathbb{P}^1$ in \cite{cmr}, linking the appropriate classical and tropical moduli spaces of curves.  
 
The arguments of Section \ref{integrals} feature wall crossing structures and scattering diagrams, generalizations of those found in \cite{MSP2}.  It seems clear that there should be many other similar enhancements, and a system for classifying these may help to uncover some sort of limiting enumerative scattering structure.    Related frameworks have now been explored in depth in \cite{cps} \cite{GHK} \cite{GHKK} and elsewhere; it would be interesting to explore the relationship of these works to this paper.

The change of coordinates given by the mirror map induces a substitution that recovers a larger sector of the big phase space than is usually explicitly encoded in Givental's $J$-function.  It may be possible to derive further relationships from the type of formalism used here.  

Finally, the combinatorial flavor of these results prompts one to seek a deeper explanation.  It would be especially gratifying to reverse our reduction to the tropical setting and explain them from a classical standpoint.

\subsection{Acknowledgements}
 I would like to thank my Ph.D. advisor, Mark Gross, for his many helpful suggestions,  reading of a draft of this paper, and support during my thesis work, from which some of this paper is drawn.  Furthermore, the exposition and techniques herein are heavily indebted to \cite{kan} and especially \cite{MSP2}. Portions of this work were completed at the Fields Institute, University of Alberta, and KU Leuven which provided excellent working conditions.  Thanks are also due to Angela Hicks and Emily Leven for their suggestions on the proof of Claim \ref{angela}.

\section{Preliminary definitions}\label{def}
Set $M:=\ZZ^2,$ $M_\RR :=M\otimes_\ZZ \RR$, $N:=\Hom(M,\mathbb{Z})$, and $\langle\cdot  ,\cdot \rangle:N\times M\rightarrow \mathbb{Z}$ the usual pairing.  Let $\Sigma$ be the toric fan of $X_\Sigma:=\mathbb{P}^2$ in $M_\mathbb{R}$.  Explicitly, let $\mathbf{m}_0:=(-1,-1)$, $\mathbf{m}_1:=(1,0),$ $\mathbf{m}_2:=(0,1)\in M$, $\rho_i=\{x\in M_\mathbb{R}|x=r\mathbf{m}_i\text{ for some }r\geq 0\}$, and $\sigma_{i,j}:=\{\mathbf{m}\in M_\mathbb{R}|\mathbf{m}=a\mathbf{m}_i+b\mathbf{m}_j\text{ for some }a, b\geq 0\}$. Then $\Sigma$ is the set of rational convex cones in $M_\mathbb{R}$ given by
$$\Sigma=\left\{\{0\}, \rho_0, \rho_1, \rho_2, \sigma_{0,1}, \sigma_{1,2}, \sigma_{2,0}\right\}.$$
We stratify $\Sigma$ by dimension, defining $\Sigma^{[0]}:=\left\{ \{0\} \right\}$, $\Sigma^{[1]}:=\left\{ \rho_0, \rho_1, \rho_2 \right\}$, and $\Sigma^{[2]}:=\left\{ \sigma_{0,1}, \sigma_{1,2}, \sigma_{2,0}\right \}$.  There is a natural filtration
$$S_i:=\bigcup_{j=0}^i \Sigma^{[j]}.$$

Denote by $T_\Sigma$ the free abelian group generated by $\Sigma^{[1]}$ and $T_\Sigma^+\subseteq T_\Sigma$ its associated semigroup.   For $\rho_i\in \Sigma^{[1]}$, denote by $t_i$ the corresponding generator in $T_\Sigma$.  

There is a natural map $p$ taking an element in $T_\Sigma$ to the corresponding linear combination of primitive generators in $M$.  Define $p:T_\Sigma\rightarrow M$ by $$t_{i }\mapsto \mathbf{m}_i.$$  Define $\overline{z}=\sum_{\rho_i\in \Sigma^{[1]}}t_i,$ and for $z=\sum_{\rho_i\in \Sigma^{[1]}} a_i t_i \in T_\Sigma$, let $|z|:=\sum_{\rho_i \in \Sigma^{[1]}} a_i \in \mathbb{Z}$.

We set $k\in \mathbb{Z}_{>0}$, which will serve as an ``order of approximation" and allow us to avoid issues of infinity in our tropical structures.  An ordered set of points $A:=\{Q, P_1, P_2, \ldots P_k\}\subset M_\mathbb{R}$ will be called an \emph{arrangement}.  For an arrangement $A$ and $Q'\in M_\mathbb{R}$, denote by $A(Q')$ the arrangement formed by replacing $Q\in A$ by $Q'$. We will often need a notion of generality, which depends on context.  In this paper, generality will always (in a fairly obvious way) refer to conditions defined by the complements of finite sets of tropical curve-like objects.  We leave it to readers to satisfy themselves with the details.  For an arrangement $A$, define $S_i(A)$ to be the translation of $S_i$ centered at $Q\in A$.

Let $\mathcal{R}_k:=\prod_{i=1}^k \{0, 1, \ldots, k\}$.  For a \emph{vector} $r=(r_1, r_2,\ldots, r_k )\in \mathcal{R}_k$, denote by $r_i$ the $i$-th entry and $\#(r)$ the number of non-zero entries of $r$.  Furthermore, let  $r\{i\}$ indicate the position and $r(i)=r_{r\{i\}}$ the value of $i$-th non-zero entry in $r$ for $1\leq i \leq \#(r)$.  Define $|r|:=\sum_{j=1}^k r_i = \sum_{i=1}^{\#(r)} r(i)$.  We occasionally need component-wise operations for $r,s \in \mathcal{R}_k$: $rs:=(r_1s_1, r_2s_2, \ldots, r_ks_k)$ and $r+s:=(r_1+s_1, r_2+s_2, \ldots, r_k+s_k)$.  Let $0\in \mathcal{R}_k$ be the additive identity.  We say $r,s\in \mathcal{R}_k$ are disjoint if $rs=0$ and $r\leq s$ if, for all $1\leq i\leq k$, $r_i\leq s_i$.  Furthermore, we say $s$ dominates $r$ (written $r\prec s$) if $r\leq s$ and $r_i >0$ if $s_i >0$ for all $1\leq i \leq k$.   If  $r\leq s$, we define $s-r \in \mathcal{R}_k$ by $(s-r)_i=s_i-r_i$ if $r_i\neq 0$ and $0$ otherwise.  Set ${|r| \choose r}:={ |r| \choose r_1, \ldots, r_k}$.  For $1\leq i \leq k$, let $e^i$ denote the elementary vector with an $i$-th entry of 1 and 0 elsewhere.

We will also need an index set containing three distinct types of elements: $$\mathcal{I}:=\{x, p_1, p_2, \ldots, q_1, q_2 , \ldots\}.$$  $\mathcal{I}$ will be used to label the three types of marked points we will encounter in our construction.

The \emph{fundamental class}, \emph{point mapping}, \emph{divisor}, and \emph{dilaton axioms} of GW theory will be used frequently.  For definitions, see \cite{c+k} or \cite{kan}, Section 2.1.

We will denote by $T_i$ a positive generator of $H^{2i}(\mathbb{P}^2,\mathbb{Z})$.  
\section{Tropical objects}\label{trop}
\subsection{Definitions}
We give a slight variation on the definitions given in \cite{MSP2}, as later considerations will require  more structure.  A \emph{metric graph} is a topological realization of a graph with possible non-compact edges, and a coordinate function (homeomorphism onto its image) $l_E:E\rightarrow \mathbb{R}_{\geq 0}$ for each edge $E$, with $l_E$ surjective when $E$ is non-compact.   We will call a finite (here referring to the number of edges and vertices), connected genus-0 metric graph a \emph{frame}.  For a frame $\Gamma$, let $\Gamma^{[1]}$ be the set of edges,  $\Gamma^{[1]}_\infty$ the set of non-compact edges, $\Gamma^{[0]}$ the set of vertices, and $\Gamma^{[0]}_i$ the set of $i$-valent vertices.  
\subsubsection{Curves}
 Let $\Gamma$ be a frame for which $\Gamma^{[0]}_1=\Gamma^{[0]}_2=\varnothing$.  Assign a weight function $w: \Gamma^{[1]}\rightarrow \mathbb{Z}_{\geq 0}$ such that $w(\Gamma^{[1]}_\infty)\subseteq \{0,1\}$ and $w^{-1}(0)\subseteq \Gamma^{[1]}_\infty$, defining a weighted frame $(\Gamma, w)$.  A marking will be a bijection $marks$ from a subset $H\subseteq \mathcal{I}$ of the form $\{x, p_1, \ldots p_n, q_1, \ldots , q_m\}$ or $\{p_1, \ldots p_n, q_1, \ldots , q_m\}$ to $w^{-1}(0)$.  We will write $marks(t)\in \Gamma_\infty ^{[1]}$ as $E_t$ for $t\in H$.  
 The data $(\Gamma, w,  marks, H)$ constitutes a \emph{marked, weighted frame}.  We will suppress the dependence on the map $marks$, simply writing $(\Gamma, w, \{x,n,m\})$ when $H=\{x, p_1, \ldots p_n, q_1, \ldots , q_m\}$ and $(\Gamma, w, \{n,m\})$ when $H=\{p_1, \ldots p_n, q_1, \ldots , q_m\}.$ 
 
A \emph{parametrized tropical curve} $(\Gamma, w, h, \{x,n,m\})$ is a marked, weighted frame $(\Gamma, w, \{x,n,m\})$ and a continuous map $h:\Gamma \rightarrow M_\mathbb{R},$ smooth on the interior of each edge of weight greater than $0$, satisfying:
\begin{itemize}
\item At any point on the interior of a given edge $E$, $h_*(\partial_x)= w(E)\mathbf{v}_E$, where $x$ is the coordinate given by $l_E$ and $\mathbf{v}_E$ is a primitive vector in $M$.
\item (Balancing condition) Let $V\in\Gamma^{[0]}$, and $E_1, \ldots E_j$ be the edges adjacent to $V$.  Let $\mathbf{m}_{E_i}=\pm \mathbf{v}_{E_i} \in M$ be a primitive vector pointing away from $h(V)$ along the direction of $h(E_i)$.  Then
$$\sum_{i=1}^mw(E_i)\mathbf{m}_{E_i}=0.$$
 \end{itemize}

A \emph{tropical curve} is an equivalence class of {parametrized tropical curves} where  $\mathcal{C} =(\Gamma, w, h, \{x,n,m\})$ is equivalent to $\mathcal{C'} =(\Gamma', w', h', \{x',n',m'\})$ if there exists an isometry $\phi:\Gamma\rightarrow \Gamma'$ respecting the marking and weight data, smooth on the interior of each edge, and with $\phi \circ h'=h$. A tropical curve $\mathcal{C} =[(\Gamma, w, h, \{x,n,m\})]$ is \emph{in} $X_\Sigma$ if, for each unmarked $E\in \Gamma^{[1]}_\infty$, $h(E)$ is a translation of some $\rho_i\in \Sigma^{[1]}$.  In this case we can define its degree as
$$\Delta(\mathcal{C}):=\sum_{\rho_i\in \Sigma^{[1]}} d_i t_i\in T_\Sigma^+$$
where $d_i$ is the number of unbounded edges of $\Gamma$ that are mapped to translations of $\rho_i$ by $h$.  

The \emph{combinatorial type} of a tropical curve $\mathcal{C} =[(\Gamma, w, h, \{x, n, m\})]$ is defined as the homeomorphism class of $\Gamma$, the markings, weights, and the data $\mathbf{m}_E$ for each edge $E$.  Note that the combinatorial type and metric structure of the underlying frame determine the image of a tropical curve up to translation in $M_\mathbb{R}$.

\subsubsection{Disks}
Our strategy for counting these curves involves a similar object, the \emph{tropical disk} (modified from the definition in \cite{kan}).  A tropical disk (or simply \emph{disk}) $\mathcal{D} =[(\Gamma, w, h, \{n,m\})]$ is defined by the same collection of data as a tropical curve, where the underlying frame $\Gamma$ has precisely one univalent vertex, $V_{out}$.  The (unique) edge of $\Gamma$ attached to $V_{out}$ will be called $E_{out}$.  We do not impose the balancing condition at $V_{out}$, but do at every other vertex.  Note that $x\in \mathcal{I}$ will not be used to mark any edge. As previously mentioned, disks should be thought of as pieces of tropical curves that have been broken apart at a vertex; the point of attachment becomes $V_{out}$.  

Define $\mathbf{m}(\mathcal{D}):=w(E_{out})\mathbf{m}^{prim}(\mathcal{D})=-p(\Delta  (\mathcal{D}))\in M_\mathbb{R}$, where $\mathbf{m}^{prim}(\mathcal{D})$ is a primitive vector tangent to $h(E_{out})$ pointing away from $h(V_{out})$.    
The formalism we introduced to treat tropical curves can be extended to disks. 

\subsection{Collections}

\subsubsection{Curves}
\begin{definition}
Let $A$ be an arrangement, $\Delta\in T_\Sigma^+$, $S\subseteq M_\mathbb{R}$, $m, \nu \in \mathbb{Z}_{\geq 0}$ and $r\in \mathcal{R}_k$.  
Then we define $\mathcal{M}_{\Delta}^{\rm curve}(A, r, T_{0,tr}^m,\psi^{\nu} S)$
to be the moduli space of tropical curves $$\mathcal{C} =[(\Gamma, w, h, \{x, p_1, \ldots, p_{\#(r)}, q_1, \ldots q_m\})]$$ in $X_\Sigma,$  such that
\begin{enumerate}
\item $h(E_{p_j})=P_{r\{j\}}$.
\item If $E_x$ shares a vertex $V_l$ with $E_{p_l}$ for $1\leq l\leq \#(r)$, then
$$Val(V_l)=2+r(l)+\nu$$
and the valence of the vertex $V_j$  attached to $E_{p_j}$ for $j\in \{1,\ldots,\#(r)\}\setminus \{k\}$ is given by
$$Val(V_j)=2+r(j)$$
\item Otherwise, the valence of the vertex $V_x$ attached to $E_x$ is given by 
$Val(V_x)=\nu+3$ and $Val(V_j)=2+r(j)$ for $1\leq j \leq \#(r)$.   
\item $h(E_x)\in S$.
\item $\Delta(\mathcal{C})=\Delta$
\end{enumerate}
\end{definition}

\begin{lemma} 
\label{lem1} 
Let $r\in \mathcal{R}_k,$  $\Delta\in T_\Sigma^+$ and $A$ be a general arrangement.  For $0\leq l \leq 2$ and $p(\Delta)=0$, $\mathcal{M}_{\Delta}^{\rm curve}(A, r, T_{0,tr}^m,\psi^{\nu} S_l)$ is a polyhedral complex of dimension $|\Delta|+m-\nu-|r|+l-2$.  By the generality of the points $P_i\in A$, the same result holds if we replace $S_l$ with $S_l(A)$.
\begin{proof}
This follows from the argument of Lemma 5.11 in \cite{kan}, changing the number of bounded edges to be $|\Delta|+m+\#(r)+1-(3+\nu+\sum_{j=1}^{\#(r)} [r(j)-1])$.
\end{proof}
\end{lemma}

\subsubsection{Disks}

\begin{definition}
\label{disktree}
Let $A$ be an arrangement, $m \in \mathbb{Z}_{\geq 0}$ and $r\in \mathcal{R}_k$.  
Then we define $Disk(A, r, T_{0,tr}^m)$
to be the set of tropical disks $$\mathcal{D} =[(\Gamma, w, h, \{p_1, \ldots, p_{\#(r)}, q_1, \ldots q_m\})]$$ (note that we do not mark disks with $x\in \mathcal{I}$)  in $X_\Sigma,$  such that
\begin{enumerate}
\item $h(E_{p_j})=P_{r\{j\}}$.
\item  The valence of each vertex $V$ is:
\begin{itemize}
\item $1$ if $V=V_{out}$
\item $3$ if $V\neq V_{out}$ is not attached to $E_{p_i}$ for any $p_i$
\item $2+r(j)$ if $V\neq V_{out}$ is attached to $p_j$
\end{itemize}
\end{enumerate}
\end{definition}
\begin{definition}
Define $RootDisk(A, r, T_{0,tr}^{m})\subseteq Disk(A, r, T_{0,tr}^m)$ to be the subset of disks with $h(V_{out})=Q$.   We define $Disk(A, T_{0,tr}^{m})$ to be the union over all $r\in\mathcal{R}_k$ of the sets $Disk(A, r, T_{0,tr}^m)$, with related subset $RootDisk(A, T_{0,tr}^{m})$.
\end{definition}

\begin{definition}
Let $\mathcal{D}=[(\Gamma, w, h, \{p_1, \ldots, p_{\#(r)}, q_1, \ldots q_m\})]$ be a tropical disk in  $Disk(A, r, T_{0,tr}^m)$.  Define the \emph{flexibility} of $\mathcal{D}$ as $$F(\mathcal{D}):=|\Delta(\mathcal{D})|+m - |r|.$$
\end{definition}
\begin{lemma}
If $A$ is a general arrangement,  the set of disks $\mathcal{D}$ in $RootDisk(A, T_{0,tr}^{m})$ with $F(\mathcal{D})=n$ is an $n-1$ dimensional polyhedral complex.  The set of such disks in $Disk(A, T_{0,tr}^m)$ is an $n+1$ dimensional polyhedral complex.  
\begin{proof}
This follows from the argument of Lemma 5.6 in \cite{kan}, replacing the idea of Maslov index with flexibility and adjusting the number of bounded edges as above.

\end{proof}
\end{lemma}

\begin{definition}
Let $\mathcal{D}$ be a tropical disk in $Disk(A, r, T_{0,tr}^m)$.  We say $\mathcal{D}$ is semirigid if $F(\mathcal{D})=1$ and rigid if $F(\mathcal{D})=0$.  Note, as one degenerate example, the single semirigid disk $\mathcal{D}\in RootDisk(A, 0, T_{0,tr}^1)$.  
\end{definition}
\section{Multiplicity}

\subsection{Disks}

We will have slightly different definitions of multiplicity for semirigid and rigid disks, closely related to the famous multiplicity due to Mikhalkin \cite{Mik}.  Our approach was inspired by the methods of  \cite{M+R}.
\begin{definition}
\label{multdisks}

Let $A$ be a general arrangement and $\mathcal{D}$ a semirigid  tropical disk in $RootDisk(A, r, T_{0,tr}^{m})$.  Then $\mathcal{D}$ can be considered as a point on the interior of a moduli space $\mathcal{M}_\mathcal{D}$ of tropical disks in $X_\Sigma$ of the same combinatorial type with $h(V_{out})=Q$, whose coordinates are given by the lengths of the bounded edges $E\in \Gamma$.  Define $ev(\mathcal{D}):\mathcal{M}_{\mathcal{D}} \rightarrow M_\mathbb{R}^{\#(r)}$ by $$ev(h)=\left(h(p_1), \ldots, h(p_{\#(r)})\right).$$
For each vertex $V\in \mathcal{V}^\Gamma$, define $n_i(V)$ to be the number of unbounded rays radiating from $V$ in the direction $\mathbf{m}_i$.  Define $$lab(\mathcal{D}):= \prod_{V\in \mathcal{V}^\Gamma}\frac{1}{n_0(V)!n_1(V)!n_2(V)!}$$and $$Mult(h):=|{\rm det}(ev)|lab(\mathcal{D}),$$
where ${\rm det}(ev)$ is the determinant of the linear part of $ev$ and we set $|{\rm det}(ev)|:=1$ if $|\#(r)|=0$.

\end{definition}

\begin{definition}
Let $\mathcal{D}$ be a rigid tropical disk in $Disk(A, r, T_{0,tr}^{m})$, with $\#(r)$ necessarily greater than 1.  We modify the definition above by placing $\mathcal{D}$ into a moduli space $\mathcal{M}_{\mathcal{D}_{rigid}}$ of tropical disks sharing the same combinatorial type, length of $E_{out}$, and image $h(E_{p_1})\in M_\mathbb{R}.$  The lengths of the rest of the bounded edges give a set of coordinates.  We define $ev'(\mathcal{D}):\mathcal{M}_{\mathcal{D}_{rigid}} \rightarrow M_\mathbb{R}^{\#(r)-1}$ by $$ev'(h)=\left(h(p_2), \ldots, h(p_{\#(r)})\right).$$   As before
$$Mult(\mathcal{D}):=|{\rm det}(ev')|lab(\mathcal{D}),$$
\end{definition}
where ${\rm det}(ev')$ is the determinant of the linear part of $ev'$ we set $|{\rm det}(ev')|:=1$ if $|\#(r)|=1$.
\subsection{Curves}

\begin{definition}
Let $S\subseteq M_\mathbb{R}$ and
$$\mathcal{C} =[(\Gamma, w, h,  \{x, p_1,\ldots, p_{\#(r)}, q_1, \ldots, q_m\})]\in \mathcal{M}_{\Delta}^{\rm curve}(A, r, T_{0,tr}^m,\psi^{\nu} S).$$ 
Denote by $\Gamma_1,\ldots, \Gamma_{w}$ the closures of each of the connected components of $\Gamma\setminus E_x$, with $h_i$ being the restriction of $h$ to $\Gamma_i$.

Each disk $\mathcal{D}_i$ defined by $\Gamma_i$ and $h_i$ is viewed as being marked by those points $p\in \{p_1,\ldots, p_{\#(r)}\}$ and $q\in \{q_1, \ldots, q_m\}$ whose corresponding edges belong to $\Gamma_i$.  That is, $\mathcal{D}_i\in Disk(A, s^i, T_{0,tr}^{m_i})$ where $m_i$ counts the number of marked points $q_j$ in $\Gamma_i$ and $r^i\in \mathcal{R}_k$ is the vector of values of $r$ corresponding to the marked points $p_j$ in $\Gamma_i$.   Note that the $\sum_i m_i=m$, $r^j\in \mathcal{R}_k$ are pairwise disjoint and $\sum_{j} r^j=r$.  Denote by 
$$\overline{Dec}(\mathcal{C}):=\{\mathcal{D}_1, \ldots, \mathcal{D}_w\}$$
 the \emph{decomposition} of $\mathcal{C}$, define ${Dec}(\mathcal{C})\subset \overline{Dec}(\mathcal{C})$ to be the subset of disks that do not consist of a single marked edge, and $simpDec(\mathcal{C})\subset Dec(\mathcal{C})$ to be the subset of disks consisting of a single unmarked, unbounded edge.
\end{definition}
The following lemma spells out the structural relationship between curves and disks.

\begin{lemma}
\label{lem2}
Let $S\subseteq M_\mathbb{R}$ be a subset.  Let 
$$\mathcal{C} =[(\Gamma, w, h, \{x, p_1, \ldots, p_{\#(r)}, q_1, \ldots q_m\})]\in \mathcal{M}_{\Delta}^{\rm curve}(A, r, T_{0,tr}^m,\psi^{\nu} S).$$ 

\begin{enumerate}
\item If $S=M_{\mathbb{R}}$ and $|r|=|\Delta|-\nu+m$, then either:
\begin{enumerate}
\item $E_x$ does not share a vertex with any $E_{p_i}$.  In this case, all but two of the disks $\mathcal{D}\in {Dec}(\mathcal{C})$ are semirigid, and the remaining two are rigid.
\item $E_x$ shares a vertex with $E_{p_j}$.  In this case, $\mathcal{D}_i\in {Dec}(\mathcal{C})$ is semirigid for all choices of $i$.
\end{enumerate}
\item If $S=C$, a general translation of $S_1$, and $|r|=|\Delta|-\nu+m-1$, then all but one of the disks $\mathcal{D}\in {Dec}(\mathcal{C})$ are semirigid, and the remaining one is rigid.
\item If $S=Q'$, a general point in $M_\mathbb{R}$, and $|r|=|\Delta|-\nu+m-2$, all disks $\mathcal{D}\in {Dec}(\mathcal{C})$ are semirigid.
\end{enumerate}
\begin{proof}
This follows from the argument of Lemma 5.12 in \cite{kan}, adjusting the dimensional requirements as dictated by Lemma \ref{lem1}.
\end{proof}
\end{lemma}

The following, rather mysterious, multiplicities are taken from \cite{MSP2} and are necessary for defining our descendent tropical invariants.  Let $\mathcal{C}$ be a tropical curve, with vertex $V_x$ attached to $E_x$.  
Define:
\begin{align*}
Mult^0_x(\mathcal{C})=&\frac{1}{n_0(V_x)!n_1(V_x)!n_2(V_x)!}\\
Mult^1_x(\mathcal{C})=&-\frac{\sum_{j=1}^{n_0(V_x)}\frac{1}{j}+\sum_{j=1}^{n_1(V_x)}\frac{1}{j}\sum_{j=1}^{n_2(V_x)}\frac{1}{j}}{n_0(V_x)!n_1(V_x)!n_2(V_x)!}\\
Mult^2_x(\mathcal{C})=&\frac{\left(\sum_{l=0}^2\sum_{j=1}^{n_l(V_x)}\frac{1}{j}\right)^2+\sum_{l=0}^2\sum_{j=1}^{n_l(V_x)}\frac{1}{j^2}}{2n_0(V_x)!n_1(V_x)!n_2(V_x)!}\\
\end{align*}

 \begin{definition}\label{tropmult}
Fix a general arrangement $A=\{Q, P_1, \ldots, P_k\}$. Let $r\in \mathcal{R}_k$, $n=\#(r)$, and $a_i:=r(i)-1$.  Recall the definition $\overline{z}=t_{0}+t_{1}+t_{2}\in T_\Sigma$.  We now define tropical curve counts that we will call \emph{descendent tropical invariants}, though they are not \emph{a priori} invariant or related to classical GW theory.  These properties will be shown in later sections.

\begin{enumerate}
\item
When $3d-2-\nu+m-|r|=0$, we define 
$$\langle \psi^{a_1}P_{r\{1\}}\ldots,\psi^{a_n} P_{r\{n\}}, T_{0,tr}^m,\psi^\nu S_0(A)\rangle_{0,d}^{trop}$$
to be 
$$\sum_\mathcal{C} Mult(\mathcal{C})$$
where the sum is over all $\mathcal{C}\in \mathcal{M}_{d\overline{z}}^{\rm curve}(A, r, T_{0,tr}^m,\psi^{\nu} S_0(A))$ with
$$Mult(\mathcal{C}):=Mult^0_x(\mathcal{C})\prod_{\mathcal{D}_i\in {Dec}(\mathcal{C})}Mult(\mathcal{D}_i).$$
\item
When $3d-1-\nu+m-|r|=0$, we define 
$$\langle \psi^{a_1}P_{r\{1\}},\ldots,\psi^{a_n} P_{r\{n\}}, T_{0,tr}^m,\psi^\nu S_1(A)\rangle_{0,d}^{trop}$$
to be 
$$\sum_\mathcal{C} Mult(\mathcal{C})$$
where the sum is over all marked tropical rational curves satisfying one of the following conditions:
\begin{enumerate}
\item $\nu\geq 0,$ 
$$\mathcal{C}\in \mathcal{M}_{d\overline{z}}^{\rm curve}(A, r, T_{0,tr}^m,\psi^{\nu} S_1(A)),$$ and no $\mathcal{D}\in {simpDec{(\mathcal{C})}}$ maps into the connected component of $S_1(A)\setminus\{Q\}$ containing $h(E_x)$.  By Lemma \ref{lem2}, there is precisely one rigid $\hat{\mathcal{D}}\in {Dec}(\mathcal{C})$.  Suppose that the connected component of $S_1(A)\setminus\{Q\}$ is $Q+\mathbb{R}_{\geq 0}\mathbf{m}_i$.  Then we define:

$$Mult(\mathcal{C}):=|\mathbf{m}(\hat{\mathcal{D}})\wedge \mathbf{m}_i|Mult^0_x(\mathcal{C})\prod_{\mathcal{D}\in {Dec}(\mathcal{C})}Mult(\mathcal{D}).$$

\item
$\nu\geq 1$ and 
$$\mathcal{C}\in \mathcal{M}_{d\overline{z}}^{\rm curve}(A, r, T_{0,tr}^m,\psi^{\nu-1} S_0(A))$$
 In this case, 
$$Mult(\mathcal{C}):=Mult^1_x(\mathcal{C})\prod_{\mathcal{D}\in {Dec}(\mathcal{C})}Mult(\mathcal{D})$$
\end{enumerate}
\item
When $3d-\nu+m-|r|=0$, we define 
$$\langle \psi^{a_1}P_{r\{1\}},\ldots,\psi^{a_n} P_{r\{n\}}, T_{0,tr}^m, \psi^\nu S_2(A) \rangle_{0,d}^{trop}$$
to be 
$$\sum_\mathcal{C} Mult(\mathcal{C})$$
where the sum is over all marked tropical curves $\mathcal{C}$ satisfying one of the following conditions:
\begin{enumerate}
\item $\nu\geq 0,$
$$\mathcal{C}\in \mathcal{M}_{d\overline{z}}^{\rm curve}(A, r, T_{0,tr}^m,\psi^{\nu} S_2(A))$$
and $E_x$ does not share a vertex with any of the $E_{p_i}$'s.  Furthermore, no $\mathcal{D}\in simpDec(\mathcal{C})$ maps into the connected component of $M_\mathbb{R}\setminus S_2(A)$ containing $h(E_x)$.  By Lemma \ref{lem2}, there are precisely two rigid disks in $Dec(\mathcal{C})$, which we call $\mathcal{D}_1$ and $\mathcal{D}_2$.   Then

$$Mult(\mathcal{C}):=|\mathbf{m}(\mathcal{D}_1)\wedge \mathbf{m}(\mathcal{D}_2)|Mult^0_x(\mathcal{C})\prod_{\mathcal{D}\in {Dec}(\mathcal{C})}Mult(\mathcal{D}).$$ 
\item \label{new} $\nu \geq 0,$
$$\mathcal{C}\in \mathcal{M}_{d\overline{z}}^{\rm curve}(A, r, T_{0,tr}^m,\psi^{\nu} S_2(A))$$
and $E_x$ shares a vertex with $E_{p_i}$.   Suppose $l$ elements of $simpDec({\mathcal{C}})$ map into the connected component of $M_\mathbb{R}\setminus L$ containing $h(E_x)$.  Then we define:

$$Mult(\mathcal{C}):={a_i+\nu-l \choose \nu}Mult^0_x(\mathcal{C})\prod_{\mathcal{D}\in {Dec}(\mathcal{C})}Mult(\mathcal{D}).$$
\item
$\nu\geq 1$ and 
$$\mathcal{C}\in \mathcal{M}_{d\overline{z}}^{\rm curve}(A, r, T_{0,tr}^m,\psi^{\nu-1} S_1(A))$$
Furthermore, no $\mathcal{D}\in simpDec(\mathcal{C})$ maps into the connected component of $S_1(A)\setminus\{Q\}$ containing $h(E_x)$.  By Lemma \ref{lem2}, there is precisely one rigid $\hat{\mathcal{D}}\in Dec(\mathcal{C})$.  Suppose that the connected component of $S_1(A)\setminus\{Q\}$ is $Q+\mathbb{R}_{\geq 0}\mathbf{m}_i$.  Then we define:

$$Mult(\mathcal{C}):=|\mathbf{m}(\hat{\mathcal{D}})\wedge \mathbf{m}_i|Mult^1_x(\mathcal{C})\prod_{\mathcal{D}\in {Dec}(\mathcal{C})}Mult(\mathcal{D}).$$

\item $\nu\geq 2$ and 
 $$\mathcal{C} \in \mathcal{M}_{d\overline{z}}^{\rm curve}(A, r, T_{0,tr}^m,\psi^{\nu-2} S_0(A))$$
 In this case, 
 $$Mult(\mathcal{C}):=Mult^2_x(\mathcal{C})\prod_{\mathcal{D}\in {Dec}(\mathcal{C})}Mult(\mathcal{D}).$$
 \end{enumerate}
 \end{enumerate}
 If $3d-\nu+m-|r|+i\neq 2$ (we will call this \emph{incompatible dimension}), we define $$\langle \psi^{a_1}P_{r\{1\}},\ldots,\psi^{a_n} P_{r\{n\}}, T_{0,tr}^m, \psi^\nu S_i(A) \rangle_{0,d}^{trop}:=0$$
\end{definition}
\begin{definition} 
 For $\sigma\in \Sigma$, define
 $$\langle \psi^{a_1}P_{r\{1\}},\ldots,\psi^{a_n} P_{r\{n\}}, T_{0,tr}^m,\psi^\nu S_i(A)\rangle_{d,\sigma}^{trop}$$
to be the contribution to $\langle \psi^{a_1}P_{r\{1\}},\ldots,\psi^{a_n} P_{r\{n\}}, T_{0,tr}^m,\psi^\nu S_i(A)\rangle_{0,d}^{trop}$
from curves with $h(E_x)$ mapping to the interior of $Q+\sigma$.
\end{definition}

We make a simple observation regarding these curve counts.
\begin{lemma}
\label{tropfun}
The descendent tropical invariants described above satisfy a tropical fundamental class axiom:
\begin{align*}
\langle \psi^{a_1}P_{r\{1\}},\ldots&,\psi^{a_n} P_{r\{n\}}, T_{0,tr}^m,\psi^\nu S_i(A)\rangle_{0,d}^{trop}\\ 
=&\sum_{j=1}^n \langle \psi^{a_1}P_{r\{1\}},\ldots, \psi^{a_j-1}P_{r\{j\}}, \ldots,\psi^{a_n} P_{r\{n\}}, T_{0,tr}^{m-1},\psi^\nu S_i(A)\rangle _{0,d}^{trop}\\ 
&+ \langle \psi^{a_1}P_{r\{1\}},\ldots, \psi^{a_n} P_{r\{n\}}, T_{0,tr}^{m-1},\psi^{\nu-1} S_i(A)\rangle_{0,d}^{trop},
\end{align*}
where the above counts are taken to be zero if any of the exponents on $\psi$ are negative.
\begin{proof}
This is immediate if any (and thus all) of the counts appearing are of incompatible dimension.  Otherwise, this can be seen by removing the edge $E_{q_m}$ from each of the curves contributing to the invariant on the left hand side, thereby generating curves contributing to invariants appearing on the right hand side.  In most cases, the multiplicity remains unchanged, and so the LHS is easily described in terms of the RHS.  In the curves appearing in part \ref{new} of the above definition, equality of contributions follows from the familiar identity
$${a+1\choose b}={a\choose b}+{a\choose b-1}.$$  See Figure \ref{ex1} for an example.
\end{proof}
\end{lemma}
 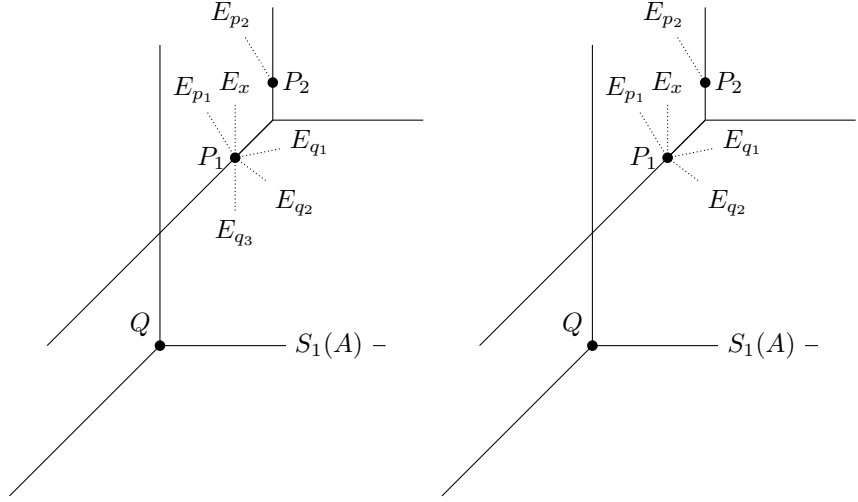
\begin{figure}[ht]
 \centering
 \subfloat{\begin{tikzpicture}
    \coordinate (Q) at (0,0) ;
    \coordinate (P_1) at (1,2.5) ;
    \coordinate (P_2) at (1.5,3.5);
    \coordinate (v_0) at (-1,-1);
     \coordinate (v_1) at (1,0);
     \coordinate (v_2) at (0,1);
      \coordinate (v_3) at (1,1);
      \coordinate (int) at (1.5, 3);
\coordinate (p) at (-.55,.9);
\coordinate (x) at (0,1);
\coordinate (q1) at (1, .2);
\coordinate(q2) at (.8, -.6);
\coordinate(q3) at (0,-1);

    \draw (0,0) -- (3,0) node [near end, fill=white] {$S_1(A)$};
    \draw (0,0)-- (-2,-2) ;
    \draw (0,0) -- (0,4) ;
\draw (P_2) -- ($(P_2)+(v_2)$);
\draw (P_1) -- (int);
\draw (P_2) -- (int);
\draw (int) -- ($(int)+2*(v_1)$);
\draw (int) -- ($(int)+3*(v_0)$);
\draw[densely dotted] (P_1) -- ($(P_1)+(p)$) node [at end, fill=white] {$E_{p_1}$};
\draw[densely dotted] (P_2) -- ($(P_2)+(p)$) node [at end, fill=white] {$E_{p_2}$};
\draw[densely dotted] (P_1) -- ($(P_1)+(x)$) node [at end, fill=white] {$E_{x}$};
\draw[densely dotted] (P_1) -- ($(P_1)+(q1)$) node [at end, fill=white] {$E_{q_1}$};
\draw[densely dotted] (P_1) -- ($(P_1)+(q2)$) node [at end, fill=white] {$E_{q_2}$};
\draw[densely dotted] (P_1) -- ($(P_1)+(q3)$) node [at end, fill=white] {$E_{q_3}$};

\fill[black] (Q) circle (2pt) node[above left] {$Q$};
\fill[black] (P_1) circle (2pt) node[left] {$P_1$};
\fill[black] (P_2) circle (2pt) node[right] {$P_2$};
 \end{tikzpicture} 
 } 
\subfloat{ \begin{tikzpicture}
    \coordinate (Q) at (0,0) ;
    \coordinate (P_1) at (1,2.5) ;
    \coordinate (P_2) at (1.5,3.5);
    \coordinate (v_0) at (-1,-1);
     \coordinate (v_1) at (1,0);
     \coordinate (v_2) at (0,1);
      \coordinate (v_3) at (1,1);
      \coordinate (int) at (1.5, 3);
\coordinate (p) at (-.55,.9);
\coordinate (x) at (0,1);
\coordinate (q1) at (1, .2);
\coordinate(q2) at (.8, -.6);
\coordinate(q3) at (0,-1);

    \draw (0,0) -- (3,0) node [near end, fill=white] {$S_1(A)$};
    \draw (0,0)-- (-2,-2) ;
    \draw (0,0) -- (0,4) ;
\draw (P_2) -- ($(P_2)+(v_2)$);
\draw (P_1) -- (int);
\draw (P_2) -- (int);
\draw (int) -- ($(int)+2*(v_1)$);
\draw (int) -- ($(int)+3*(v_0)$);
\draw[densely dotted] (P_1) -- ($(P_1)+(p)$) node [at end, fill=white] {$E_{p_1}$};
\draw[densely dotted] (P_2) -- ($(P_2)+(p)$) node [at end, fill=white] {$E_{p_2}$};
\draw[densely dotted] (P_1) -- ($(P_1)+(x)$) node [at end, fill=white] {$E_{x}$};
\draw[densely dotted] (P_1) -- ($(P_1)+(q1)$) node [at end, fill=white] {$E_{q_1}$};
\draw[densely dotted] (P_1) -- ($(P_1)+(q2)$) node [at end, fill=white] {$E_{q_2}$};

\fill[black] (Q) circle (2pt) node[above left] {$Q$};
\fill[black] (P_1) circle (2pt) node[left] {$P_1$};
\fill[black] (P_2) circle (2pt) node[right] {$P_2$};
 \end{tikzpicture}  }
\caption{
 Two tropical curves.  The dotted lines indicate unbounded edges of weight zero.  By Definition \ref{tropmult} \ref{new}, the curve on the left contributes ${5\choose 2}$ to $\langle \psi^3 P_1, P_2, T_{0,tr}^3,\psi^2 M_\mathbb{R}\rangle_{0,1}^{trop},$ while the curve on the right contributes ${4 \choose 2}$ to $\langle \psi^2 P_1, P_2, T_{0,tr}^2,\psi^2 M_\mathbb{R}\rangle_{0,1}^{trop}$ and ${4\choose 1}$ to $\langle \psi^3 P_1, P_2, T_{0,trop}^2,\psi^1 M_\mathbb{R}\rangle_{0,1}^{trop}$.}
 \label{ex1}
 \end{figure}
 
\section{Tropical B-model} \label{bmod}
\subsection{Tropical Landau-Ginzburg potential}
Let $A$ be a general arrangement.  We give a framework, generalized from that appearing in \cite{MSP2}, yielding a more refined Landau-Ginzburg potential whose integral recovers tropical versions of a broader class of GW invariants.
\begin{definition}
To $P_i\in A$ associate the variables $u_{i,j}$ in the ring:
$$R_{k}:=\frac{\mathbb{C}[\{u_{i,j}\}_{i,j}]}{I}$$
with $j\in \mathbb{Z}_{\geq 0}$ and $i \in\mathbb{Z}_{> 0}$, where $I$ is the ideal generated by the set
$$\{u_{i,j}u_{i,j'}|1\leq i \leq k, \text{ }0\leq j \leq j' \leq k\}\cup \{u_{i,j}| i>k \text{ or } j\geq k\}$$
Let $\overline{m}\in \mathbb{Z}_{\geq 0}$ and define
$$\mathfrak{R}_{k,\overline{m}}:=R_k[y_{0,0}]/(y_{0,0}^{\overline{m}+1}).$$
\end{definition}
For $r\in \mathcal{R}_k$, define
$$u_r:=\prod_{i=1}^{\#(r)} u_{r\{i\},r(i)-1}.$$
Let
$$y_{2,j}:=\sum_{i=1}^k u_{i,j}.$$
Note that $a\in \mathfrak{R}_{k,\overline{m}}$ can be uniquely represented as 
$$a=\sum_{\substack{0\leq m \leq \overline{m}\\r\in \mathcal{R}_k}}a_{r,m} u_r y_{0,0}^m$$
with $a_{r,m}\in \mathbb{C}$, where we abuse notation by writing $u_ry_{0,0}^m$ to denote its equivalence class in $\mathfrak{R}_{k,\overline{m}}$.
\begin{definition}
Let $\mathcal{D}$ be a  tropical disk in $Disk(A, r, T_{0,tr}^{w})$.  Define $u_{\mathcal{D}}:=u_r$ and $y_{0,0}^\mathcal{D} := \frac{y_{0,0}^w}{w!}.$

\end{definition}
\begin{definition}
Let $\mathcal{D}$ be a rigid or semirigid disk.  Then
$$Mono(\mathcal{D}):=Mult(\mathcal{D})u_{\mathcal{D}}z^{\Delta(\mathcal{D})}y_{0,0}^\mathcal{D}\in \mathbb{C}[T_\Sigma]\otimes_\mathbb{C}\mathfrak{R}_{k,\overline{m}}$$
where $z^{\Delta(\mathcal{D})}\in\mathbb{C}[T_{\Sigma}]$ is the monomial associated to $\Delta(\mathcal{D})\in T_{\Sigma}$.  We will write $x_i=z^{t_{i}}$, so $z^{n_0{t_{0}}+n_1{t_{1}}+n_2{t_{2}}} =x_0^{n_0}x_1^{n_1}x_2^{n_2}$.
\end{definition}
\begin{definition}
We define the $(k,\overline{m})$ descendent Landau Ginzburg potential associated to $A$ as 
$$W_{k,\overline{m}}(A):=\sum_\mathcal{D} Mono(\mathcal{D})$$
where the sum is over all semirigid disks in $RootDisk(A, r, T_{0,tr}^{m})$ for any $r \in \mathcal{R}_k$ and $m\leq \overline{m}$.
\end{definition}

\begin{definition}
$$W_{basic}(A):=x_0+x_1+x_2$$
\end{definition}

\subsection{B-model moduli}
Here we review the construction given in \cite{MSP2}.   Recall the  map $p:T_\Sigma\rightarrow M$ given by $p(t_{i})=\mathbf{m}_i$.  As $p$ is surjective (a consequence of non-singularity), we have the following exact sequence:
$$0\rightarrow K_\Sigma\rightarrow T_\Sigma \xrightarrow p M\rightarrow 0$$
with $K_\Sigma$ the kernel of $p$.  
Dualizing over $\mathbb{Z}$ gives
$$0\rightarrow \Hom(M, \ZZ) \rightarrow {\rm Hom}_\mathbb{Z}(T_\Sigma, \mathbb{Z})\rightarrow {\rm Pic} X_\Sigma \rightarrow 0$$
Tensoring with $\mathbb{C}^\times$ gives the sequence

 $$0\rightarrow \Hom(M, \ZZ)\otimes\mathbb{C}^\times \rightarrow {\rm Hom}(T_\Sigma, \mathbb{C}^\times)\xrightarrow\kappa {\rm Pic} X_\Sigma\otimes\mathbb{C}^\times \rightarrow 0$$
 defining $\kappa$, which provides the family of mirrors to $X_\Sigma$.  
 Set $$\check{\mathcal{X}}:=\Hom(T_\Sigma, \mathbb{C}^\times)=\Spec\mathbb{C}[T_\Sigma].$$

The \emph{K\"ahler moduli space} of $X_\Sigma$ is defined to be 
$$ \mathcal{M}_\Sigma:={\rm Pic}X_\Sigma\otimes \mathbb{C}^\times=\Spec\mathbb{C}[K_\Sigma].$$
In our case, $K_\Sigma\cong\mathbb{Z}$.  Note that $\kappa$, by definition, is now a map:
$$\kappa:\Spec\mathbb{C}[T_\Sigma]\rightarrow \mathcal{M}_\Sigma$$
A fiber of $\kappa$ over a closed point of $\mathcal{M}_\Sigma$ is isomorphic to $\Spec\mathbb{C}[M]$.  
Define the $(k,\overline{m})$-order thickening of the K\"ahler moduli space by
$$\mathcal{M}_{k}:=\mathcal{M}_\Sigma \times \Spec \mathfrak{R}_{k,\overline{m}}$$
and likewise
$$\check{\mathcal{X}}_{\Sigma, k}:=\check{\mathcal{X}}_\Sigma \times \mathfrak{R}_{k,\overline{m}}.$$
This yields a family
$$\kappa:\check{\mathcal{X}}_{\Sigma, k}\rightarrow \mathcal{M}_{\Sigma, k}$$
By construction, $W_{k,\overline{m}}(A)$ is a regular function on $\check{X}_{\Sigma, k}$, and should be considered as a family of Landau-Ginzburg potentials.

\section{Integrals}\label{integrals}
In this section, we will give the main result of \cite{thesis} and a summary of the methods used in its proof.  Elements of the argument which can be easily generalized from those found in  \cite{kan} are given with a reference to the relevant result, while subtler points are presented in more detail.
Define
$$\gamma_{a,tr}:=\sum_{1 \leq v+1,w \leq k}\psi^v P_w u_{w,v}$$ as a formal expression for insertion into tropical invariants, to be expanded linearly.  For example, 
$$\langle \gamma_{a,tr}, S_2(A)\rangle^{trop}_{0,d}:=\sum_{1\leq j+1, i\leq k}\langle  \psi^j P_i, S_2(A)\rangle^{trop}_{0,d} u_{i,j}.$$

\begin{theorem} \label{theorem1}
A choice of a general arrangement $A$ gives rise to a function $W_{k,0}(A) \in \mathbb{C}[T_\Sigma]\otimes_\mathbb{C} \mathfrak{R}_{k,\overline{m}}$, and hence a family of Landau-Ginzburg potentials on the family $\check{\mathcal{X}}_{\Sigma, k}\xrightarrow \kappa \mathcal{M}_{\Sigma, k}$ with a relative nowhere-vanishing two-form $\Omega$.  This data gives rise to a local system $\mathcal{R}$ on $\mathcal{M}_{\Sigma, k}\otimes \Spec\mathbb{C}[\hbar, \hbar^{-1}]$, whose fiber over $(\kappa, \hbar)$ is $H_2((\check{\mathcal{X}}_{\Sigma, k})_\kappa, \rm{Re}(W_{basic}(A)/\hbar)\ll0)$.  Letting $y_{1,0}:=\log (\kappa)$, there exists a multi valued basis $\Xi_0, \Xi_1, \Xi_2$ of $\mathcal{R}$ satisfying the requirements of Section 1 of \cite{MSP2} such that
\begin{align*}
\sum_{i=0}^2 \alpha^i \int_{\Xi_i} e^{W_{k,0}(A)/\hbar} \Omega=\hbar^{-3\alpha} \sum_{j=0}^2  \left(\alpha\hbar\right)^j e^{y_{1,0}\alpha} \Theta_j
\end{align*}
where we have identified a fiber of $\mathcal{R}^\vee$ with $\mathbb{C}[\alpha]/(\alpha^3)$, with $\alpha^i$ dual to $\Xi_i$.  Then

\begin{align*}
\Theta_0: &=1+\sum_{d>0,w\geq 0} \frac{\hbar^{-1}}{w!}\langle \frac{S_0(A)}{\hbar - \psi}, \gamma_{a,tr}^w \rangle^{trop}_{0,d}e^{y_{1,0}d}\\
\Theta_1: &=\sum_{d>0,w\geq 0} \frac{\hbar^{-1}}{w!}\langle \frac{S_1(A)}{\hbar - \psi}, \gamma_{a,tr}^w \rangle^{trop}_{0,d}e^{y_{1,0}d}\\
\Theta_2: &=\hbar^{-1}\sum_{j=0}^k y_{2,j}(-\hbar)^j+\sum_{d>0,w\geq 0} \frac{\hbar^{-1}}{w!}\langle \frac{S_2(A)}{\hbar - \psi}, \gamma_{a,tr}^w \rangle^{trop}_{0,d}e^{y_{1,0}d}.
\end{align*}
Furthermore, the result does not depend on the choice of $A$.
\begin{proof}
See \cite{thesis} and below.
\end{proof}
\end{theorem}

\subsection{Scattering diagrams}
The first step in the proof of Theorem \ref{theorem1} is to construct a set of structures that govern the combinatorics of the potential $W_{k,0}(A)$. These methods are part of a larger theory developed by Kontsevich, Soibelman, Gross, Siebert, and a number of collaborators; although it will not be apparent here, there are deep and unexpected links to other areas of mathematics (see \cite{GHKK}).  The incarnation we use is particularly simple.   

One can form an object $\mathcal{T}$ called a \emph{tropical tree} from a rigid tropical disk  $\mathcal{D} =[(\Gamma, w, h, \{p_1, \ldots, p_{\#(r)}, q_1, \ldots q_m\})]$ in $Disk(A, r, T_{0,tr}^m)$ by deleting the vertex $V_{out}$ from the underlying frame (thereby creating a non-compact edge $E_{out}$ with $w(E_{out})\in \mathbb{Z}_{>0}$) and modifying $h$ by extending the image of $E_{out}$ to be an unbounded ray in $M_{\mathbb{R}}$.  We denote by $Tree(A,r,T_{0,tr}^m)$ the set of such trees derived from rigid disks in $Disk(A,r, T_{0,tr}^m)$. Note that the subset of rigid disks in $Disk(A, r, T_{0,tr}^m)$ has a natural fibration over a finite set $Tree(A,r,T_{0,tr}^m)$. We associate the same multiplicity $Mult(\mathcal{T})$, monomial $Mono(\mathcal{T})$, and flexibility $F(\mathcal{T})=0$ to a tree as we do to any of its overlying disks.  There is a finite set of trees associated to a general arrangement:
$$\mathfrak{T}(A)_{k,\overline{m}}=\bigcup_{w\leq \overline{m}} Tree(A,r,T_{0,tr}^m).$$
If one represents the set $\mathfrak{T}(A)_{k,\overline{m}}$ in $M_\mathbb{R}$ by drawing the outgoing edge corresponding to each rigid tree, a striking pattern emerges.  The points at which these outgoing edges intersect have rays sprouting from them, as rigid trees can be glued at such a point to form a ``child" tree.  The weight and direction of the outgoing edge of the  child is, by the balancing condition, determined by the weights and directions of its parents' outgoing edges.  Similarly, the multiplicity and monomial of the child tree are simply determined by those of its parents.
This ``scattering" at points of intersection gives our tool its name.  We hereafter specialize to the case of $\overline{m}=0$, which was addressed in \cite{thesis}, although it is largely straightforward to generalize to $\overline{m}>0$.
\begin{definition}
The following definition is from \cite{MSP2}.
\begin{enumerate}
\item A \emph{ray} or \emph{line} is a pair $(\mathfrak{d}, f_\mathfrak{d})$ such that
	\begin{itemize}
		\item $\mathfrak{d} \subseteq M_\mathbb{R}$ is given by
		$$\mathfrak{d} = \mathbf{m}_{init}-\mathbb{R}_{\geq 0}p(m_0)$$
		if $\mathfrak{d}$ is a ray and
		$$\mathfrak{d} = \mathbf{m}_{init}-\mathbb{R}p(m_0)$$
		if $\mathfrak{d}$ is a line, where $\mathbf{m}_{init} \in M_\mathbb{R}$ with $m_0\in T_\Sigma$ satisfying $$-\mathbf{m}_\mathfrak{d}:=p(m_0)\neq 0.$$  The set $\mathfrak{d}$ is the \emph{support} of the ray or line.  If $\mathfrak{d}$ is a ray, then $\mathbf{m}_{init}$ is called the \emph{initial point} and is denoted $Init(\mathfrak{d})$.
		\item $f_\mathfrak{d}\in \mathbb{C}[z^{m_0}]\otimes_\mathbb{C} \mathfrak{R}_{k,\overline{m}}$
		\item $f_\mathfrak{d}\equiv 1 \mod (\{u_{i,j}\}_{i,j})z^{m_0}$
	\end{itemize}
\item A \emph{scattering digram} $\mathfrak{D}$ is a finite collection of lines and rays.
\end{enumerate}
We will sometimes write $w(\mathfrak{d}):=w(E_{out})$ for walls $\mathfrak{d}$ in $\mathfrak{D}(A)_{k,0}$.

If $\mathfrak{D}$ is a scattering diagram, we write
$$Supp(\mathfrak{D}):=\bigcup_{\mathfrak{d}\in \mathfrak{D}} \mathfrak{d} \subseteq M_\mathbb{R}$$
and
$$Sing(\mathfrak{D}):=\bigcup_{\mathfrak{d}\in \mathfrak{D}} \partial \mathfrak{d} \cup \bigcup_{\substack{\mathfrak{d}_1, \mathfrak{d}_2\\ {\rm dim} \mathfrak{d}_1\cap \mathfrak{d}_2=0}} \mathfrak{d}_1\cap \mathfrak{d}_2$$
where $\partial\mathfrak{d}=\{Init(\mathfrak{d})\}$ if $\mathfrak{d}$ is a ray, and is empty if it is a line.
\end{definition}

\begin{definition}
We build our diagram $\mathfrak{D}(A)_{k,0}$ from the outgoing edges of the trees in $\mathfrak{T}(A)_{k,0}$.  The ray in $\mathfrak{D}(A)_{k,0}$ corresponding to a tree $\mathcal{T}\in \mathfrak{T}(A)_{k,0}$ is of the form $(\mathfrak{d}, f_\mathfrak{d})$, where 
\begin{itemize}
\item $\mathfrak{d}=h(E_{out})$
\item $f_\mathfrak{d}=1+w(E_{out})Mult(\mathcal{T})z^{\Delta(\mathcal{T})}u_{\mathcal{T}}$
\end{itemize}
\end{definition}

\begin{definition}
Given a scattering diagram $\mathfrak{D}$ and smooth immersion $\xi:[0,1]\rightarrow M_\mathbb{R}\setminus Sing(\mathfrak{D})$ whose endpoints are not in $Supp(\mathfrak{D})$, with $\xi$ intersecting $Supp(\mathfrak{D})$ transversally, we can use this information to define a ring automorphism $\theta_{\xi,\mathfrak{D}}$ of $\mathcal{R}_{k,0}$.
Find numbers
$$0<s_1\leq s_2\leq \ldots \leq s_n <1$$
and elements $\mathfrak{d}_i$ such that $\xi(s_i)\in \mathfrak{d}_i$, $\mathfrak{d}_i\neq \mathfrak{d}_j$ if $i\neq j $ and $n$ is taken to be as large as possible to account for all walls of $\mathfrak{D}$ that are crossed by $\xi$.  For each $i\in \{1,\ldots, n\}$, define $\theta_{\xi,\mathfrak{d}_i}$ to be the automorphism with action
\begin{align*}
\theta_{\xi, \mathfrak{d}_i}(z^w)&= z^wf_{\mathfrak{d}_i}^{\langle \mathbf{n}_0, p(w) \rangle}\\
\theta_{\xi, \mathfrak{d}_i}(a)&= a
\end{align*}
for $w\in T_\Sigma$, $a\in \mathfrak{R}_{k,0}$, where $\mathbf{n}_0\in N$ is chosen to be primitive, annihilating the tangent space to $\mathfrak{d}_i$, and satisfying
$$\langle \mathbf{n}_0, \xi'(s_i)\rangle <0$$
Then $\theta_{\xi, \mathfrak{D}}:=\theta_{\xi, \mathfrak{d}_n}\circ\cdots \circ \theta_{\xi, \mathfrak{d}_1}$, where composition is taken from right to left.

\end{definition}
The reproductive process associated to $\mathfrak{D}(A)_{k,0}$ gives rise to a useful property that distinguishes it from scattering diagrams encountered in other contexts \cite{tropvert}.
\begin{lemma}
\label{lem51}
If $P\in Sing(\mathfrak{D}(A)_{k,0})$ is a singular point with $P\notin A$ and $\xi_p$ is a small loop around $P$, then $$\theta_{\xi_p,\mathfrak{D}(A)_{k,0}}=Id.$$
\begin{proof}
See \cite{kan}, Proposition 5.28.
\end{proof}
\end{lemma}

These automorphisms have another nice property: membership in $\mathbb{V}_{\Sigma, k}$, a group of automorphisms of $\mathbb{C}[T_\Sigma]\otimes_{\mathbb{C}} \mathcal{R}_{k,0}$ originally defined in \cite{koso06} as a set of Hamiltonian symplectomorphisms (see \cite{kan}, 5.4.2).  Significantly for us, these automorphisms preserve the choice of $\Omega$ referenced in Theorem \ref{theorem1} and, when acting on $W_{k,0}(A)$, leave the period integral unchanged.
\begin{lemma}
\label{lem52}
Let $\sigma \in \mathbb{V}_{\Sigma,k}$, $(w,h)\in \mathcal{M}_{\Sigma, k}\times\mathbb{C}^{\times}$ and suppose that $f$ is in the ideal generated by $(\{u_{i,j}\})$ in $\mathbb{C}[T_\Sigma]\otimes_\mathbb{C}\mathcal{R}_{k,0}$.  Then, for any cycle
$$\Xi \in H_2(\kappa^{-1}(w), {\rm Re}(W_{basic}/\hbar)\ll0, \mathbb{C}),$$
we have
$$\int_\Xi e^\frac{(W_{basic}+f)}{\hbar}\Omega=\int_\Xi e^\frac{\theta(W_{basic}+f)}{\hbar}\Omega.$$
\begin{proof}
See \cite{kan} Lemma 5.40.
\end{proof}
\end{lemma}
\subsubsection{Broken lines}
The technique of \emph{broken lines} connects $\mathfrak{D}_{k,0}(A)$ to the potential $W_{k,0}(A)$.  Every semirigid disk in $RootDisk(A, T_{0,tr}^0)$ can be uniquely described as a central (infinitely long) stem onto which a number of rigid disks are grafted.  At each point of grafting, the stem bends in a way dictated by the balancing condition.  This is easily understood in terms of the scattering diagram, because the possible points at which any particular tree can be attached (as a rigid disk) to a stem are given by the wall it contributes to $\mathfrak{D}(A)_{k,0}$.  Therefore, in order to understand semirigid disks contributing to $W_{k,0}(A)$, it is sufficient to analyze the behavior of these stems (\emph{broken lines}) with respect to the scattering diagram.  We use the following definition, adapted from \cite{MSP2}.
\begin{definition}
A \emph{broken line} with basepoint $Q'\in M_\mathbb{R}$ is a continuous proper piecewise linear map
$$\beta:(-\infty, 0] \rightarrow M_\mathbb{R}$$
with endpoint $Q'=\beta(0)$, along with some additional data.  Let
$$-\infty = s_0 <s_1<\cdots  <s_n =0$$
be the smallest set of real numbers such that $\beta|_{(s_{i-1}, s_i)}$ is linear.  Then, for each $1\leq i\leq n$, we are given the additional data of a monomial $c_iz^{w_i^\beta}\in \mathbb{C}[T_\Sigma]\otimes_\mathbb{C}\mathfrak{R}_{k,0}$ with $w_i^\beta \in T_\Sigma\setminus K_\Sigma$, satisfying:
\begin{enumerate}
\item For each $i$, $p(w_i^\beta)=-\beta'(s)$ for $s\in (s_{i-1}, s_i)$.
\item $w_1^\beta = t_i$ for some $0\leq i \leq 2 $ and $c_1=1$.
\item $\beta(s_i)\in Supp(\mathfrak{D}(A)_{k,0})\setminus Sing(\mathfrak{D}(A)_{k,0})$ for $1\leq i \leq n$.
\item If $\beta(s_i) \in \mathfrak{d}_1\cap \cdots \cap \mathfrak{d}_n$, then $c_{i+1}z^{w^\beta_{i+1}}$ is a term in 
$$(\theta_{\beta, \mathfrak{d}_1}\circ\cdots\circ\theta_{\beta, \mathfrak{d}_n})(c_iz^{w_i^\beta})$$
More explicitly, suppose that $f_{\mathfrak{d}_j}=1+c_{\mathfrak{d}_j}z^{w_{\mathfrak{d}_j}}$, $1\leq j \leq n$, with $c^2_{\mathfrak{d}_j}=0$, and $\mathbf{n}\in N$ is primitive, orthogonal to all of the $\mathfrak{d}_j$'s, and chosen so that
\begin{align*}
(\theta_{\beta, \mathfrak{d}_1}\circ\cdots\circ\theta_{\beta, \mathfrak{d}_1})(c_iz^{w_i^\beta}) &= c_iz^{w_i^\beta}\prod_{j=1}^n (1+c_{\mathfrak{d}_j}z^{w_{\mathfrak{d}_j}})^{\langle \mathbf{n}, p(w_i^\beta)\rangle}\\
&= c_iz^{w_i^\beta}\prod_{j=1}^n (1+{\langle \mathbf{n}, p(w_i^\beta)\rangle c_{\mathfrak{d}_j}z^{m_{\mathfrak{d}_j}}}).
\end{align*}
Then we must have 
$$c_{i+1}z^{w^\beta_{i+1}}=\prod_{j\in J}{\langle \mathbf{n}, p(w_i^\beta)\rangle c_{\mathfrak{d}_j}z^{m_{\mathfrak{d}_j}}}$$
for some $J\subseteq \{1,\ldots, n\}$.  We interpret this as $\beta$ being bent at time $s_i$ by $\mathfrak{d}_j$ for $j\in J$.
\end{enumerate}

\end{definition}
\begin{proposition} If $A$ is a general arrangement, there is a one-to-one correspondence between broken lines with endpoint $Q$ and semirigid disks in $RootDisk(A, T_{0,tr}^0).$  In addition, if $\beta$ is a broken line corresponding to a disk $\mathcal{D}$, and $cz^w$ is the monomial associated to the last segment of $\beta$, then
$$cz^w=Mono(\mathcal{D})$$
\begin{proof}
See Proposition 5.32 of \cite{kan}.
\end{proof}
\end{proposition}
\subsection{Wall crossing and evaluation of integrals}
To evaluate the integral appearing in Theorem \ref{theorem1}, we must first show that changing the arrangement $A$ transforms $W_{k,0}(A)$ by the action of an element of $\mathbb{V}_{\Sigma, k}$, and will thus leave the integral unchanged.   We examine the effect on the integral by replacing $A$ by $A(Q')$ while moving $Q'$ out to infinity in a particular direction, noting then that the contribution to the integral from terms with certain monomials vanishes.  We can then understand the contribution of these monomials to the integral associated to $A$ by considering the wall crossing automorphisms we encounter as we move $Q'$ back to $Q$.  These automorphisms, and thus the period integrals, can be interpreted in terms of tropical curves.  Using this technique, Theorem \ref{theorem1} was proven by Gross in the non-descendent case (in our notation, $r=(r_1, \ldots, r_n)$  with $r_i\leq 1$) in \cite{kan}.  The same techniques are modified to treat the descendent case (arbitrary $r_i$) in \cite{thesis}.  This modification is straightforward in most cases, as the relevant scattering diagrams have identical structure away from the points in the arrangement $A$.  
\subsection{Wall crossing}

\begin{lemma}\label{diskstotree}

Let $Q'\in M_\mathbb{R}$ be very near $P_l$ and let $\mathcal{D}_j\in RootDisk(A(Q'), r^j, T_{0,tr}^{0})$ for $1\leq j \leq n$ be semirigid  disks such that the vectors $r^j$ are pairwise disjoint and $e^lr^j=0$ for all $j$ .  If $\sum_{j=1}^{n}\mathbf{m}(\mathcal{D}_j)\neq 0$, then the disks $\mathcal{D}_j$ can be joined at $P_l$ to give a rigid tree $\mathcal{T} \in Tree(A,r,T_{0,tr}^0)$ with outgoing edge $P_l+\mathbb{R}\sum_{j=1}^{n}\mathbf{m}(\mathcal{D}_j)$, where $r:=ne^l+\sum_{j=1}^n r^j$.  Let $M_i\subseteq\{\mathcal{D}_j\}$ be the set of our original disks which are simply outgoing edges in the direction $\mathbf{m}_i$.  Then
$$Mult(\mathcal{T})=\frac{u_{k,n-1}}{|M_0|!|M_1|!|M_2|!}\prod_{i\in\{1,\ldots, n\}}Mult(\mathcal{D}_i)$$
\begin{proof}
It is easy to see that the resulting tree is rigid.  The rest follows from linear algebra.
\end{proof}
\end{lemma}
\begin{lemma}\label{treetodisks}
Let $\mathcal{T}\in Tree(A,r,T_{0,tr}^0)$ with $r_l=n+1$.  Then, by splitting $\mathcal{T}$ at the vertex $V$ mapping to $P_l$, we can form $n$ semirigid tropical disks rooted at some $Q'\in M_\mathbb{R}$, chosen near $P_l$.
\begin{proof}
Call the $n$ tropical disks formed by the above procedure $\mathcal{D}_1, \ldots, \mathcal{D}_n$, with $\mathcal{D}_j\in RootDisk(A(Q'), r^j, T_{0,tr}^0)$.  Each $F(\mathcal{D}_j)\leq 1$ as $\mathcal{T}$ is rigid.  Note $F(\mathcal{T})=|\Delta(\mathcal{T})|-|r|=0$, $|r|=\nu+\sum_{j=1}^\nu |r^j|$ and $|\Delta(\mathcal{T})|=\sum_{j=1}^\nu |\Delta(\mathcal{D}_j)|$, so 
$$
\sum_{j=1}^n F(\mathcal{D}_j) =n.
$$

Thus $F(\mathcal{D}_j)=1$ for all $j\in\{1,\ldots,n\}$.
\end{proof}
\end{lemma}

\begin{figure}[ht]

 \centering
 \subfloat{\begin{tikzpicture}
    \coordinate (P_w) at (0,0) ;
    \coordinate (P_1) at (1,2.5) ;
    \coordinate (P_2) at (1.5,3.5);
    \coordinate (v_0) at (1,1);
    \coordinate (v_1) at (1,1);
     \coordinate (v_2) at (-1,1);
     \coordinate (v_3) at (-1,0);
      \coordinate (v_4) at (-1,-1);
     \coordinate (v_5) at (1,-1);
     \coordinate (disk0) at (1,0);
      \coordinate (out) at (0,-1); 
           
 \fill[black] (P_w) circle (2pt) node[above] {$P_l$};

\draw(P_w) -- ($(P_w)+2*(v_1)$) node [at end, fill=white] {$\mathcal{D}_1$};
\draw(P_w) -- ($(P_w)+2*(v_2)$) node [at end, fill=white] {$\mathcal{D}_2$};
\draw(P_w) -- ($(P_w)+2*(v_3)$) node [at end, fill=white] {$\mathcal{D}_3$};
\draw(P_w) -- ($(P_w)+2*(v_4)$) node [at end, fill=white] {$\mathcal{D}_4$};
\draw(P_w) -- ($(P_w)+2*(v_5)$) node [at end, fill=white] {$\mathcal{D}_5$};
\draw[dir](P_w) -- ($(P_w)+2*(disk0)$) node [at end, fill=white] {$\mathfrak{d}_0$};
\draw[densely dotted]($(disk0)+(v_1)$)--(disk0) node[at start, fill= white]{$B_0$};
\draw[dir, densely dotted](disk0)--($(disk0)+2*(out)$);

 \end{tikzpicture} 
 }
 \subfloat{\begin{tikzpicture}
    \coordinate (P_w) at (0,0) ;
    \coordinate (P_1) at (1,2.5) ;
    \coordinate (P_2) at (1.5,3.5);
    \coordinate (v_0) at (1,1);
    \coordinate (v_1) at (1,1);
     \coordinate (v_2) at (-1,1);
     \coordinate (v_3) at (-1,0);
      \coordinate (v_4) at (-1,-1);
     \coordinate (v_5) at (1,-1);
     \coordinate (disk0) at (1,0);
     \coordinate (disk4) at (-.5,-1);
          \coordinate (disk5) at (.5,-1);
      \coordinate (out) at (0,-1); 

 \fill[black] (P_w) circle (2pt) node[above] {$P_l$};
\draw(P_w) -- ($(P_w)+2*(v_1)$) node [at end, fill=white] {$\mathcal{D}_0, \mathcal{D}_1$};
\draw(P_w) -- ($(P_w)+2*(v_2)$) node [at end, fill=white] {$\mathcal{D}_2$};
\draw(P_w) -- ($(P_w)+2*(v_3)$) node [at end, fill=white] {$\mathcal{D}_3$};
\draw(P_w) -- ($(P_w)+2*(v_4)$) node [at end, fill=white] {$\mathcal{D}_4$};
\draw[dir](P_w) -- ($(P_w)+2*(disk5)$) node [at end, fill=white] {$\mathfrak{d}_5$};
\draw[densely dotted]($.5*(disk5)+(v_5)$)--($.5*(disk5)$) node[at start, fill= white]{$B_5$};
\draw[dir, densely dotted]($.5*(disk5)$)--($.5*(disk5)+1.5*(out)$);
\end{tikzpicture}}

  \subfloat{\begin{tikzpicture}
    \coordinate (P_w) at (0,0) ;
    \coordinate (P_1) at (1,2.5) ;
    \coordinate (P_2) at (1.5,3.5);
    \coordinate (v_0) at (1,1);
    \coordinate (v_1) at (1,1);
     \coordinate (v_2) at (-1,1);
     \coordinate (v_3) at (-1,0);
      \coordinate (v_4) at (-1,-1);
     \coordinate (v_5) at (1,-1);
     \coordinate (disk0) at (1,0);
      \coordinate (out) at (0,-1); 

 \fill[black] (P_w) circle (2pt) node[above] {$P_l$};
\draw(P_w) -- ($(P_w)+2*(v_0)$) node [at end, fill=white] {$\mathcal{D}_0$};
\draw(P_w) -- ($(P_w)+2*(v_2)$) node [at end, fill=white] {$\mathcal{D}_2$};
\draw(P_w) -- ($(P_w)+2*(v_3)$) node [at end, fill=white] {$\mathcal{D}_3$};
\draw(P_w) -- ($(P_w)+2*(v_4)$) node [at end, fill=white] {$\mathcal{D}_4$};
\draw(P_w) -- ($(P_w)+2*(v_5)$) node [at end, fill=white] {$\mathcal{D}_5$};
\draw[dir](P_w) -- ($(P_w)+2*(disk0)$) node [at end, fill=white] {$\mathfrak{d}_1$};
\draw[densely dotted]($(disk0)+(v_1)$)--(disk0) node[at start, fill= white]{$B_1$};
\draw[dir, densely dotted](disk0)--($(disk0)+2*(out)$);

 \end{tikzpicture} 
 }
 \subfloat{\begin{tikzpicture}
    \coordinate (P_l) at (0,0) ;
    \coordinate (P_1) at (1,2.5) ;
    \coordinate (P_2) at (1.5,3.5);
    \coordinate (v_0) at (1,1);
    \coordinate (v_1) at (1,1);
     \coordinate (v_2) at (-1,1);
     \coordinate (v_3) at (-1,0);
      \coordinate (v_4) at (-1,-1);
     \coordinate (v_5) at (1,-1);
     \coordinate (disk0) at (1,0);
      \coordinate (out) at (0,-1); 
     
 \fill[black] (P_l) circle (2pt) node[above] {$P_l$};
           
\draw(P_l) -- ($(P_l)+2*(v_1)$) node [at end, fill=white] {$\mathcal{D}_0, \mathcal{D}_1$};
\draw[dir](P_l) -- ($(P_l)+2*(v_3)$) node [at end, fill=white] {$\mathcal{D}_3, \mathfrak{d}_2$};
\draw(P_l) -- ($(P_l)+2*(v_4)$) node [at end, fill=white] {$\mathcal{D}_4$};
\draw(P_l) -- ($(P_l)+2*(v_5)$) node [at end, fill=white] {$\mathcal{D}_5$};
\draw[densely dotted]($(v_3)+(v_2)$)--(v_3) node[at start, fill= white]{$B_2$};
\draw[dir, densely dotted](v_3)--($(v_3)+2*(out)$);
\end{tikzpicture}}

 \subfloat{\begin{tikzpicture}
    \coordinate (P_w) at (0,0) ;
    \coordinate (P_1) at (1,2.5) ;
    \coordinate (P_2) at (1.5,3.5);
    \coordinate (v_0) at (1,1);
    \coordinate (v_1) at (1,1);
     \coordinate (v_2) at (-1,1);
     \coordinate (v_3) at (-1,0);
      \coordinate (v_4) at (-1,-1);
     \coordinate (v_5) at (1,-1);
     \coordinate (disk0) at (1,0);
      \coordinate (out) at (0,-1); 
           
 \fill[black] (P_w) circle (2pt) node[above] {$P_l$};

\draw(P_w) -- ($(P_w)+2*(v_1)$) node [at end, fill=white] {$\mathcal{D}_0, \mathcal{D}_1$};
\draw(P_w) -- ($(P_w)+2*(v_2)$) node [at end, fill=white] {$\mathcal{D}_2$};
\draw[dir](P_w) -- ($(P_w)+2*(v_4)$) node [at end, fill=white] {$\mathcal{D}_4, \mathfrak{d}_3$};
\draw(P_w) -- ($(P_w)+2*(v_5)$) node [at end, fill=white] {$\mathcal{D}_5$};
\draw[densely dotted]($.5*(v_4)+(v_3)$)--($.5*(v_4)$) node[at start, fill= white]{$B_3$};
\draw[dir, densely dotted]($.5*(v_4)$)--($.5*(v_4)+1.5*(out)$);
\end{tikzpicture}}
 \subfloat{\begin{tikzpicture}
    \coordinate (P_w) at (0,0) ;
    \coordinate (P_1) at (1,2.5) ;
    \coordinate (P_2) at (1.5,3.5);
    \coordinate (v_0) at (1,1);
    \coordinate (v_1) at (1,1);
     \coordinate (v_2) at (-1,1);
     \coordinate (v_3) at (-1,0);
      \coordinate (v_4) at (-1,-1);
     \coordinate (v_5) at (1,-1);
     \coordinate (disk0) at (1,0);
     \coordinate (disk4) at (-.5,-1);
      \coordinate (out) at (0,-1); 
           
 \fill[black] (P_w) circle (2pt) node[above] {$P_l$};

\draw(P_w) -- ($(P_w)+2*(v_1)$) node [at end, fill=white] {$\mathcal{D}_0, \mathcal{D}_1$};
\draw(P_w) -- ($(P_w)+2*(v_2)$) node [at end, fill=white] {$\mathcal{D}_2$};
\draw(P_w) -- ($(P_w)+2*(v_3)$) node [at end, fill=white] {$\mathcal{D}_3$};
\draw(P_w) -- ($(P_w)+2*(v_5)$) node [at end, fill=white] {$\mathcal{D}_5$};
\draw[dir](P_w) -- ($(P_w)+2*(disk4)$) node [at end, fill=white] {$\mathfrak{d}_4$};
\draw[densely dotted]($($.5*(disk4)$) +(v_4)$)--($.5*(disk4)$)  node[at start, fill= white]{$B_4$};
\draw[dir, densely dotted]($.5*(disk4)$) --($.5*(disk4)+1.5*(out)$);
\end{tikzpicture}}
 
\caption{
 An example of the behavior encountered in the proof of Theorem \ref{indthm}.  The first three broken lines sare on the right hand side of the wall, while the last three are on the left.}
  \label{deg}
 \end{figure}
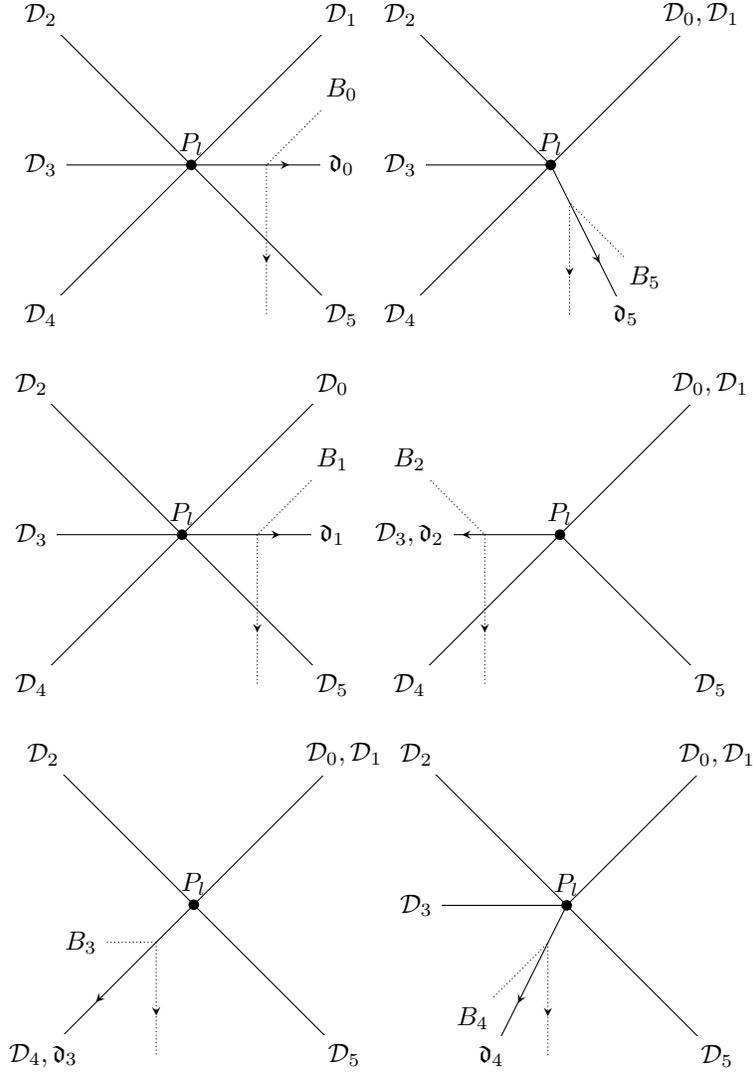

\begin{theorem}\label{indthm}
If $A(Q)$ and $A(Q')$ are two general arrangements and $\xi$ is a path connecting $Q$ and $Q'$ for which $\theta_{\xi,\mathfrak{D}(A)_{k,0}}$ is defined, $$\theta_{\xi,\mathfrak{D}(A)_{k,0}}(W_{k,0}(A(Q)))=W_{k,0}(A(Q')).$$
\begin{proof}
This theorem, except for one case, follows from a slight modification of the argument found in \cite{kan}, Theorem 5.35.  The strategy is to analyze the behavior of so-called \emph{degenerate} broken lines.  These occur as the limits of deformations of ordinary broken lines; as one deforms the base point, two bends can can converge to a single point on the broken line, or one of the bends can approach a singular point of the scattering diagram.  See Definition 5.34 of \cite{kan} for a formal definition.  One subdivides the plane by a set of walls composed of those from $\mathfrak{D}_{0,k}(A)$ in addition to those formed by such degenerate broken lines; the change in $W_{k,0}\left (A\left (\xi(s)\right) \right)$ as $\xi(s)$ crosses one of these walls will be seen to be generated by an automorphism of $\mathbb{C}[T_\Sigma]\otimes_{\mathbb{C}} \mathfrak{R}_{k,0}$.  This automorphism can be understood as a type of mutation process on the broken lines with endpoint $\xi(s)$.  The only case which requires an argument significantly different from that appearing in \cite{kan} is an analysis of the autormorphisms induced by degenerate broken lines which bend at some $P_l\in A$. 

For $\hat{Q}\in M_\mathbb{R}$,  denote by $\mathfrak{B}(\hat{Q})$ the set of broken lines in  $\mathfrak{D}(A)_{k,0}$ with endpoint $\hat{Q}$. Suppose $\xi(s_0)$ is in some wall $L$ to which $\xi$ is transverse, and for small $\epsilon>0$, let $Q_1:=\xi(s_0-\epsilon)$ and $Q_2:=\xi(s_0+\epsilon)$.  Let $\mathbf{n}\in N$ be a primitive vector annihilating the tangent space to $L$ at $\xi(s_0)$ and taking a smaller value on $Q_1$ than $Q_2$.  We decompose $\mathfrak{B}(Q_i)$ into $\mathfrak{B}^+(Q_i)$, $\mathfrak{B}^0(Q_i)$, and $\mathfrak{B}^-(Q_i)$, where the membership of $\beta \in \mathfrak{B}(Q_i)$ is determined the sign of $\langle \beta_*(-\partial/\partial s|_{s=s_0}), \mathbf{n}\rangle.$

These decompositions allows us to write $$W_{k,0}(A(Q_i))=W^-_{k,0}(A(Q_i))+W^0_{k,0}(A(Q_i))+W^+_{k,0}(A(Q_i)).$$ Following the techniques in \cite{kan}, one can show 
$$\theta_{\xi',\mathfrak{D}(A)_{k,0}}(W^\pm_{k,0}(A(Q_1)))=W^\pm_{k,0}(A(Q_2)),$$ where $\xi'$ is the segment of $\xi$ joining $Q_1$ to $Q_2$.  

For the remaining case,  we will partition $\mathfrak{B}(Q_i)^0=\bigsqcup_{j=1}^l \mathfrak{B}^i_j$ and show that for each $j\in \{1, \ldots, l\}$, $\mathfrak{B}^1_j$ and $\mathfrak{B}^2_j$ make equal contributions to $W_{k,0}(Q_1)$ and $W_{k,0}(Q_2)$ respectively.  We will assume that a broken line with endpoint $\xi(s_0)$ passes through at most one singular point.  The general case follows by an induction argument.

Suppose $\beta_1 \in \mathfrak{B}(Q_1)^0$ deforms continuously to $\beta_2\in \mathfrak{B}(Q_2)^0$.  In this case, each $\beta_i$ will appear in a one element set, say $\mathfrak{B}^i_j$, and each $\mathfrak{B}^i_j$ will make the same contribution to $W_{k,n}(Q_i)$.

If $\beta\in \mathfrak{B}(Q_1)^0$ cannot be continuously deformed to an element of $ \mathfrak{B}(Q_2)^0$, then it must deform to a degenerate broken line when the base point reaches $\xi(s_0)$.  In other words, there is a map $B:(-\infty, 0]\times [0,s_0]\rightarrow M_\mathbb{R}$ such that $B|_{(-\infty, 0]\times [0,s_0)}$ is a continuous deformation of $\beta$ and $B|_{(-\infty, 0]\times \{s_0\}}:=\beta'$ is a degenerate broken line bending at $P\in Sing( \mathfrak{D}(A)_{k,0})$ at time $s'$.  There are two cases to examine: $P\in \{P_1,\ldots, P_k\}$ and $P\notin \{P_1,\ldots, P_k\}$.  We explain the former, which requires a more sophisticated argument than that appearing in \emph{loc. cit.}

Suppose $P=P_l$ and select $\hat{Q}$ very near $P_l$.  We know that $\beta$ bends along exactly one ray $\mathfrak{d}_0$ emanating from $P_l$ whose attached function has a monomial containing $u_{l,w}$.  By construction, $\mathfrak{d}_0$ is produced by a rigid tropical tree, which, by Lemma \ref{treetodisks}, is constructed from $w+1$ semirigid descendent tropical disks with endpoint $\hat{Q}$, which we will call $\mathcal{D}_1, \ldots, \mathcal{D}_{w+1}$, with $\mathcal{D}_j\in RootDisk(A(\hat{Q}), r^j, T_{0,tr}^0)$.  Also note that $B|_{(-\infty, s']\times\{s_0\}}$ is a broken line ending at $P_l$, corresponding to a semirigid disk $\mathcal{D}_0\in RootDisk(A(\hat{Q}), r^0, T_{0,tr}^0).$  The vectors $r^j$ are disjoint for all $0\leq j \leq w+1$ .  We can expect to form something like a rigid tropical tree $\mathcal{T}_j$ for each $0\leq j \leq w+1$ by joining all of the $\mathcal{D}_i$ except for $\mathcal{D}_j$ at $P_l$ and extending an unbounded outgoing edge $\mathfrak{d}_j$ as dictated by the balancing condition.  See Figure \ref{deg}.  We may happen to have $\sum_{l\neq j} \mathbf{m}(\mathcal{D}_l)=0$ and the result will not strictly qualify as a rigid descendent tropical tree, but these exceptional cases won't be problematic.  Let $M_i\subseteq \{\mathcal{D}_0, \ldots \mathcal{D}_{w+1}\}$ be the subset of disks that are simply unbounded rays pointing in the direction $\mathbf{m}_i$ from $P_l$.  Each choice of $0 \leq i \leq w+1$ where $w(\mathfrak{d}_i)\neq 0$ gives rise to a broken line $B_i$ bending at $\mathfrak{d}_i$ constructed from the concatenation of the broken line defining $\mathcal{D}_j$ and $B|_{[s',0]\times\{s_0\}}$.  We will show that the contributions to $W_{k,0}(Q_1)$ and $W_{k,0}(Q_2)$ from associated broken lines are equal.  Notice that the side of the wall that each $B_i$ inhabits is dictated by the sign of $\mathbf{m}_{\mathfrak{d}_i}\wedge \mathbf{m}(B|_{[s',0]\times\{s_0\}})$, where $\mathbf{m}(B|_{[s',0]\times\{s_0\}})$ gives the direction vector for the outgoing piece of the broken line.  Furthermore, $\mathbf{m}(B|_{[s',0]}\times\{s_0\})$ is given by $\sum_{j=0}^{w+1}\mathbf{m}(\mathcal{D}_j)$.

The monomial obtained from the bend of $B_i$ at $\mathfrak{d}_i$ is given by
\begin{align*}
w(\mathfrak{d}_i)\langle \mathbf{n}_i, \mathbf{m}(\mathcal{D}_i)\rangle &Mono(\mathcal{D}_i) Mono(\mathcal{T}_i)= \\
&\frac{w(\mathfrak{d}_i)\langle \mathbf{n}_i, \mathbf{m}(\mathcal{D}_i)\rangle  Mono(\mathcal{T}_i)}{|M_0\setminus\{\mathcal{D}_i\}|!|M_1\setminus\{\mathcal{D}_i\}|!|M_2\setminus\{\mathcal{D}_i\}|!} \prod_{n\neq i} Mono(\mathcal{D}_n)\\
= &\frac{w(\mathfrak{d}_i)\langle \mathbf{n}_i, \mathbf{m}(\mathcal{D}_i)\rangle u_{l,w}}{|M_0\setminus\{\mathcal{D}_i\}|!|M_1\setminus\{\mathcal{D}_i\}|!|M_2\setminus\{\mathcal{D}_i\}|!} \prod_{n} Mono(\mathcal{D}_n),
\end{align*}
where $\mathbf{n}_i\in N$ is orthogonal to $\mathfrak{d}_i$ and chosen so that 
$$w(\mathfrak{d}_i)\langle \mathbf{n}_i, \mathbf{m}(\mathcal{D}_i)\rangle = |\left(\sum_{n\neq i}\mathbf{m}(\mathcal{D}_n)\right)\wedge \mathbf{m}(\mathcal{D}_i)|$$
 (as $\mathbf{m}_{\mathfrak{d}_i}$ is given by $\sum_{n\neq i}\mathbf{m}(\mathcal{D}_n)$.\\
The result then follows from some basic observations.  First, 
\begin{align*}
0=\left(\sum_{j=0}^{w+1} \mathbf{m}(\mathcal{D}_j)\right)^{\wedge 2}&=\sum_{j=0}^{w+1} \mathbf{m}(\mathcal{D}_j)\wedge\left(\sum_{n=0}^{w+1} \mathbf{m}(\mathcal{D}_n)\right)\\&=\sum_{j=0}^{w+1} \mathbf{m}(\mathcal{D}_j)\wedge\left(\sum_{n\neq j} \mathbf{m}(\mathcal{D}_n)\right).
\end{align*}
Let $I^-:=\left\{n\in \{0,\ldots, k+1\right\}| \left(\sum_{n\neq i}\mathbf{m}(\mathcal{D}_n)\right)\wedge \mathbf{m}(\mathcal{D}_i))<0\}$ under the identification of $\wedge^2 M_{\mathbb{R}}$ with $\mathbb{Z}$, with $I^0$ and $I^+$ defined analogously.  Then \\

\begin{align*}
0=&\sum_{j \in I^-} \mathbf{m}(\mathcal{D}_j)\wedge\left(\sum_{n\neq j} \mathbf{m}(\mathcal{D}_n)\right)+\sum_{j \in I^+} \mathbf{m}(\mathcal{D}_j)\wedge\left(\sum_{n\neq j} \mathbf{m}(\mathcal{D}_n)\right)\\&+\sum_{j \in I^0} \mathbf{m}(\mathcal{D}_j)\wedge\left(\sum_{n\neq j} \mathbf{m}(\mathcal{D}_n)\right)\\
=&\sum_{j \in I^-} \mathbf{m}(\mathcal{D}_j)\wedge\left(\sum_{n\neq j} \mathbf{m}(\mathcal{D}_n)\right)+\sum_{j \in I^+} \mathbf{m}(\mathcal{D}_j)\wedge\left(\sum_{n\neq j} \mathbf{m}(\mathcal{D}_n)\right).
\end{align*}

A series of implications follows:
 \begin{align*}
 -\sum_{j \in I^-} \mathbf{m}(\mathcal{D}_j)\wedge\left(\sum_{n\neq j} \mathbf{m}(\mathcal{D}_n)\right)&=\sum_{j \in I^+} \mathbf{m}(\mathcal{D}_j)\wedge\left(\sum_{n\neq j} \mathbf{m}(\mathcal{D}_n)\right)\\\sum_{j \in I^-} \left|\mathbf{m}(\mathcal{D}_j)\wedge\left(\sum_{n\neq j} \mathbf{m}(\mathcal{D}_n)\right)\right|&=\sum_{j \in I^+} \left|\mathbf{m}(\mathcal{D}_j)\wedge\left(\sum_{n\neq j} \mathbf{m}(\mathcal{D}_n)\right)\right|\\\sum_{j \in I^-} w(\mathfrak{d}_j)\langle \mathbf{n}_j, \mathbf{m}(\mathcal{D}_i)\rangle&=\sum_{j \in I^+} w(\mathfrak{d}_j)\langle \mathbf{n}_j, \mathbf{m}(\mathcal{D}_i)\rangle.
\end{align*}
Therefore
\begin{align}
\sum_{j \in I^-} \frac{w(\mathfrak{d}_j)\langle \mathbf{n}_j, \mathbf{m}(\mathcal{D}_i)\rangle u_{l,w}}{|M_0|!|M_1|!|M_2|!} \prod_{n} {\rm Mono}(\mathcal{D}_n) \label{keyeq}\\=\sum_{j \in I^+}  \frac{w(\mathfrak{d}_j)\langle \mathbf{n}_j, \mathbf{m} (\mathcal{D}_i)\rangle u_{l,w}}{|M_0|!|M_1|!|M_2|!} \prod_{n} {\rm Mono}(\mathcal{D}_n).\nonumber
\end{align}
Equation \ref{keyeq} closely resembes our dseired result, as $I^+$ indexes disks related to broken lines contributing to one of $W_{k,0}(A(Q_1))$, $W_{k,0}(A(Q_2))$ and $I^-$ indexes those which contribute to the other.  To conclude, note that at most one broken line is produced for each set $M_j$, so we can say that the contribution from each $B_i$ (where $\mathcal{D}_i\in M_j$) is just $\frac{1}{|M_j|}$ of the contribution from the unique broken line produced by $M_j$.  That is, the contribution from $B_i\in M_j$ should be considered as 
\begin{eqnarray}
\nonumber \frac{1}{|M_j|}  \frac{w(\mathfrak{d}_i)\langle \mathbf{n}_i, \mathbf{m}(\mathcal{D}_i)\rangle u_{l,w}}{|M_0\setminus\{\mathcal{D}_i\}|!|M_1\setminus\{\mathcal{D}_i\}|!|M_2\setminus\{\mathcal{D}_i\}|!} \prod_{n} {\rm Mono}(\mathcal{D}_n)\\=\frac{w(\mathfrak{d}_i)\langle \mathbf{n}_i, p(\Delta (\mathcal{D}_i)\rangle u_{l,w}}{|M_0|!|M_1|!|M_2|!} \prod_{n} {\rm Mono}(\mathcal{D}_n)\nonumber
\end{eqnarray}
Of course, if $\mathcal{D}_i\notin \cup_j M_j$ then the contribution is
\begin{eqnarray}
\nonumber \frac{w(\mathfrak{d}_i)\langle \mathbf{n}_i, \mathbf{m}(\mathcal{D}_i)\rangle u_{l,w} }{|M_0\setminus\{\mathcal{D}_i\}|!|M_1\setminus\{\mathcal{D}_i\}|!|M_2\setminus\{\mathcal{D}_i\}|!} \prod_{n} {\rm Mono}(\mathcal{D}_n)\\= \frac{w(\mathfrak{d}_i)\langle \mathbf{n}_i, \mathbf{m}(\mathcal{D}_i)\rangle u_{l,w}}{|M_0|!|M_1|!|M_2|!} \prod_{n} {\rm Mono}(\mathcal{D}_n)\nonumber
\end{eqnarray}
Thus, \ref{keyeq} shows that the sum of the monomials generated by our set of broken lines on either side of the wall is equal. Deforming any of the $B_i$ to degenerate at $P_l$ will result in the same scenario, showing that broken lines degenerating at $P_l$ (for a particular deformation of $Q$) can be partitioned into sets which give equal contributions to $W_{k,0}(A(Q_1))$ and $W_{k,0}(A(Q_2))$.  As $\theta_{\xi,\mathfrak{D}(A)_{k,0}}(W^0_{k,0}(A(Q_1)))=W^0_{k,0}(A(Q_1))$, $\theta_{\xi,\mathfrak{D}(A)_{k,0}}(W^0_{k,0}(A(Q_1)))=W^0_{k,0}(A(Q_2))$.
 \end{proof}
\end{theorem}

\begin{thm}
Let $A$ and $A'$ be two general arrangements. 
$$W_{k,0}(A')=\theta(W_{k,0}(A))$$ for some $\theta \in \mathbb{V}_{\Sigma, k}$.  
\begin{proof}
This follows from a relatively straightforward generalization of the techniques of \cite{kan}, Theorem 5.39.  For details, see \cite{thesis}.
\end{proof}
\end{thm}

\subsection{Evaluation of Integrals}
The results of the previous subsection yield the following useful observation.
\begin{lemma}
\label{lem53}
For $\Xi \in H_2(\kappa^{-1}(u), {\rm Re}(W_{basic}/\hbar)\ll0, \mathbb{C}),$ the integral
$$\int_\Xi e^\frac{W_{k,0}(A)}{\hbar}\Omega$$
is independent of the choice of general arrangement $A$.
\end{lemma}
The following gives us a numerical expression of the integrals.
\begin{lemma}
\label{lem54}
Restricting to $x_0x_1x_2=\kappa$,  we have 
\begin{align*}
\sum_{i=0}^2 \alpha^i \int_{\Xi_i}e^{(x_0+x_1+x_2)/\hbar} x_0^{n_0}x_1^{n_1}x_2^{n_2}\Omega= \hbar^{-3\alpha}\kappa^{\alpha} \sum_{i=0}^2 \psi_i(n_0, n_1, n_2) \alpha ^i
\end{align*}
where $\alpha$ and $\Xi_i$ are as defined in Theorem \ref{theorem1} and $$\psi_i(n_0, n_1, n_2) = \sum_{d=0}^\infty D_i(d,n_0, n_1, n_2)\hbar^{-(3d-n_0-n_1-n_2)}\kappa^d$$
where the terms $D_i$ are numerical quantities defined in  \cite{kan}, Lemma 5.43.  For $w=n_0t_{0}+ n_1t_{1}+n_2t_{2}\in T_{\Sigma}$, we write $D_i(d,w):=D_i(d,n_0, n_1, n_2)$.
\begin{proof}
See \cite{kan}, Lemma 5.43.
\end{proof}
\end{lemma}

\begin{definition}
\label{def36}
Fix general an arrangement $A$.  For $Q'\in M_\mathbb{R}$, let $S_{k,0}(Q')$ be the finite set of triples $(c,\nu, w)$ with $c\in \mathfrak{R}_{k,0}$, $\nu \geq 0$ an integer, and $w\in T_\Sigma$ such that:
$$e^{(W_{k,0}(A(Q'))-W_{basic}(A(Q')))/\hbar} = \sum_{(c,\nu, w)\in S_{k,0}(Q')} c\hbar^{-\nu}z^w,$$
with each term $c\hbar^{-\nu}z^w$ of the form $\hbar^{-\nu}\prod_{i=1}^\nu \rm{Mono}(\mathcal{D}_i)$ for $\mathcal{D}_1,\ldots, \mathcal{D}_\nu$ semirigid disks with endpoint $Q'$.

Then
$$L_i^d(Q'):=\sum_{(c,\nu,w)\in S_{k,0}(Q')}c\hbar^{-(3d+\nu-|w|)}D_i(d,w).$$ 
\end{definition}
\begin{lemma}
\begin{align*}
\sum_{i=0}^2 \alpha^i \int_{\Xi_i}e^{W_{k,0}(A)/\hbar}= \hbar^{-3\alpha}\kappa^{\alpha} \sum_{i=0}^2 \sum_{d\geq 0} L_i^d(Q) \kappa^d \alpha ^i
\end{align*}
\begin{proof}
Follows from definitions.
\end{proof}
\end{lemma}

\begin{definition}
\label{def37}
For each cone $\sigma\in \Sigma$, $\sigma$ is the image under $p$ of a proper face $\tilde{\sigma}$ of the cone $C:=T_\Sigma^+\otimes \mathbb{R}$.  For $d\geq 0$, define $C_d\subseteq C$ to be the cube 
$$C_d=\left\{\sum_{i=0}^2 n_it_{i}| 0\leq n_i \leq d\right\}$$
and for $\sigma\in \Sigma$
$$\tilde{\sigma}_d:=(\tilde{\sigma}+C_d)\setminus\bigcup_{\tau\subsetneq\sigma, \tau\in \Sigma} (\tilde{\tau}+C_d).$$
where $+$ denotes the Minkowski sum.
\end{definition}

\begin{definition}
\label{def38}
For $\sigma \in \Sigma$ and $Q' \in M_\mathbb{R}$, define 
$$L_{i,\sigma}^d(Q'):=\sum_{(c,\nu,w)\in S_{k,0}(Q'), \, w\in \tilde{\sigma}_d}c\hbar^{-(3d+\nu-|w|)}D_i(d,w).$$ 
\end{definition}

\begin{lemma}
\label{lem55}
$L_i^d(Q')=\sum_{\sigma\in\Sigma} L_{i,\sigma}^d(Q').$
\begin{proof}
Follows immediately from definitions.
\end{proof}
\end{lemma}

\begin{lemma}
\label{lem56}
Let $\{0\}\neq \omega\in \Sigma$, and $\mathbf{v}\in \omega$ be non-zero.  Then
$$\lim_{s\rightarrow \infty}L_{i,\omega}^d(Q+s\mathbf{v})=0.$$
\begin{proof}
See \cite{kan}, Lemma 5.51. 
\end{proof}
\end{lemma}
\begin{definition}
\label{def39}
Let $\mathfrak{D}=\mathfrak{D}(A)_{k,0}$.  Let $\mathbf{C}_1$ and $\mathbf{C}_2$ be connected components of $M_\mathbb{R}\setminus\mathfrak{D}$ with $dim(\bar{\mathbf{C}}_1\cap\bar{\mathbf{C}}_2)=1$.  Pick general points $Q_i$ in $\mathbf{C}_i$, and let $\xi$ be a general path from $Q_1$ to $Q_2$ intersecting $Supp(\mathfrak{D})$ exactly once at $\xi(s_0)$, a nonsingular point of $Supp(\mathfrak{D})$.  Let $\mathfrak{d}\in \mathfrak{D}$ contain $\xi(s_0)$, and let $\mathbf{n}_\mathfrak{d}$ be a primitive vector perpendicular to $\mathfrak{d}$ pointing toward $Q_1$.   

Suppose that $f_\mathfrak{d}=1+c_\mathfrak{d}z^{m_\mathfrak{d}}$. Select $\alpha,\tau\in \Sigma$ with $dim(\tau)=dim(\alpha)+1$ and $\alpha\subseteq \tau$.  Note that there is a unique index $j\in \{0,1,2\}$ such that $\mathbf{m}_j\in \tau$ but $\mathbf{m}_j\notin \alpha$.  Call this index $j(\alpha,\tau)$.

Define  
$$L^d_{i, \mathfrak{d}, \xi, \alpha \rightarrow \tau}:=\sum_{(c,\nu, m)}cc_\mathfrak{d}\langle \mathbf{n}_\mathfrak{d}, \mathbf{m}_{j(\alpha,\tau)} \rangle  D_i(d, m+m_{\mathfrak{d}}+t_{j(\alpha,\tau)})h^{-(\nu+3d-|m+m_\mathfrak{d}|)},$$
where we sum over all $(c,\nu, m)$ in $S_{k,0}(Q_1)$ satisfying $m+m_\mathfrak{d}\in \tilde{\alpha}_d$ but $m+m_\mathfrak{d}+t_{j(\alpha,\tau)}\in \tilde{\tau}_d$.  If $(c,\nu, m)$ satisfies these condition, then we say that $c\hbar^{-\nu}z^m$ \emph{contributes} to $L^d_{i, \mathfrak{d}, \xi, \alpha \rightarrow \tau}$.\\
Define
$$L^d_{i, \xi, \alpha \rightarrow \tau}:=\sum_{\mathfrak{d}}L^d_{i, \mathfrak{d}, \xi, \alpha \rightarrow \tau}$$
where $\mathfrak{d}$ ranges over all rays of $\mathfrak{D}$ containing $\xi(s_0)$.  
In order to define this operation for a general path $\xi$, break it up into segments of the type outlined above.
\end{definition}
\begin{lemma}
\label{lem57}
Let $\xi_j$ be the straight path joining $Q$ with $Q+s\mathbf{m}_j$ for $s\gg0$.  Let $\xi_{j,j+1}$ be the loop based at $Q$ which passes linearly from $Q$ to $Q+s\mathbf{m}_j$, takes a large circular arc to $Q+s\mathbf{m}_{j+1}$, and then proceeds linearly from $Q+s\mathbf{m}_{j+1}$ to $Q$.  Here we take $j$ modulo 3, and $\xi_{j,j+1}$ is always a counterclockwise loop.   Then 
$$L_i^d(Q)=L_{i,\{0\}}^d(Q)-\sum_{j=0}^2L^d_{i,\xi_j,\{0\}\rightarrow \rho_j}-\sum_{j=0}^2L^d_{i,\xi_{j,j+1},\rho_{j+1}\rightarrow \sigma_{j,j+1}}.$$
\begin{proof}
See Lemma 5.54. of  \cite{kan}.  
\end{proof}
\end{lemma}
\begin{definition}
If $\mathcal{C}$ is a tropical curve contributing to 
$$\langle \psi^{r(1)-1}P_{r\{1\}},\ldots, \psi^{r(\#(r))-1}P_{r\{\#(r)\}}, \psi^\nu S_i(A) \rangle^{\rm trop}_{d,0},$$ 
define $r^\mathcal{C}\in \mathcal{R}_k$ to be the corresponding vector and $u_\mathcal{C}:=u_{r^\mathcal{C}}$.
\end{definition}
\begin{lemma}
\label{lem58}
\begin{align*}
L_{i, \{0\}}^d&(Q)=\delta_{0,d}\delta_{0,i}+\\&\sum_{\substack{\nu\geq i\\ r\in \mathcal{R}_k\\|r|=3d-2+i-\nu}}\langle \psi^{r(1)-1}P_{r\{1\}},\ldots, \psi^{r(\#(r))-1}P_{r\{\#(r)\}}, \psi^\nu S_i(A) \rangle^{\rm trop}_{d,\{0\}}u_{r}h^{-(\nu+2-i)}.
\end{align*}
\begin{proof}
See Lemma 5.55 of \cite{kan}.  \end{proof}
\end{lemma}

\begin{lemma}
\label{lem59}
\begin{align*}
-&L_{i, ,\xi_j,\{0\}\rightarrow \rho_j}^d=\\
&\sum_{\substack{\nu\geq i-1\\ r\in \mathcal{R}_k\\|r|=3d-2+i-\nu}}\langle \psi^{r(1)-1}P_{r\{1\}},\ldots, \psi^{r(\#(r))-1}P_{r\{\#(r)\}}, \psi^\nu S_i(A) \rangle^{\rm trop}_{d,\rho_j}u_{r}h^{-(\nu+2-i)}
\end{align*}
\begin{proof}
See Lemma 5.56 of \cite{kan}.
\end{proof}
\end{lemma}
\begin{lemma}
\label{lem60}
For each point $P\in Sing(\mathfrak{D})$, let $\xi_P$ be a small counterclockwise loop around $P$, small enough so that it doesn't go around any other point of $Sing(\mathfrak{D})$.  Then
$$L^d_{i,\xi_{j,j+1},\rho_{j+1}\rightarrow\sigma_{j,j+1} }= \sum_{P\in Sing(\mathfrak{D})\cap(Q+\sigma_{j,j+1})} L^d_{i,\xi_P, \rho_{j+1}\rightarrow \sigma_{j,j+1}}$$
\begin{proof}
See Lemma 5.57 of \cite{kan}. 
\end{proof}
\end{lemma}

\begin{lemma}
\label{lem61}
Let $P\in Sing(\mathfrak{D})\cap (Q+\sigma_{j,j+1})$, and suppose that 
$$P\notin A.$$
Then
\begin{align*}
-L_{i,\xi_P,\rho_{j+1}\rightarrow\sigma_{j,j+1}}^d
=\sum_{\nu\geq 0}\sum_\mathcal{C}Mult(\mathcal{C})u_\mathcal{C}\hbar^{-(\nu+2-i)}
\end{align*}
where the sum is over curves $\mathcal{C}$ contributing to 
$$\langle \psi^{r(1)-1}P_{r\{1\}},\ldots, \psi^{r(\#(r))-1}P_{r\{\#(r)\}},  \psi^\nu S_i(A) \rangle^{\rm trop}_{d,\sigma_{j,j+1}}$$
 for $r\in \mathcal{R}_k$ with $|r|=3d-2+i-\nu$ and $h(E_x)=P$.
\begin{proof}
See Lemma 5.58 of \cite{kan}.
\end{proof}
\end{lemma}
The following is the only place in the evaluation of the integral that requires a significant generalization of Gross's techniques.
\begin{lemma}
\label{lem62}
Let $P\in Sing(\mathfrak{D})\cap (Q+\sigma_{j,j+1})$, and suppose that $P=P_l\in A$.
Then
\begin{align*}
-&L_{i,\xi_P,\rho_{j+1}\rightarrow\sigma_{j,j+1}}^d\\
&=\sum_{w=1}^{k} u_{l,w-1}(-\hbar)^w\delta_{d,0}\delta_{2,i}+\sum_{\nu\geq 0}\sum_\mathcal{C}Mult(\mathcal{C})u_\mathcal{C}\hbar^{-(\nu+2-i)}
\end{align*}
where the sum is over curves $\mathcal{C}$ contributing to $$\langle \psi^{r(1)-1}P_{r\{1\}},\ldots, \psi^{r(\#(r))-1}P_{r\{\#(r)\}},  \psi^\nu S_i(A) \rangle^{\rm trop}_{d,\sigma_{j,j+1}}$$ for  for $r\in \mathcal{R}_k$ with $|r|=3d-2+i-\nu$ and $h(E_x)=P_l$.
\begin{proof}
Here we assume $i=2$, and write
$$L_{P,j}=L_{2,\xi_P,\rho_{j+1}\rightarrow\sigma_{j,j+1}}^d$$

Choose a basepoint $Q'$ near $P_l$.  As discussed in Lemmas \ref{treetodisks} and \ref{diskstotree}, sets of $a+1$ semirigid disks with endpoint $Q'$ not bending near $P_l$ correspond to rays in $\mathfrak{D}$ based at $P_l$ whose monomial contains $u_{l,a}$.  More precisely, sets of semirigid disks $\{\mathcal{D}_1, \ldots, \mathcal{D}_{a+1}\}$ not bending near $P_l$ with endpoint $Q'$, $\prod_iMono(\mathcal{D}_i)\neq 0,$ and $\sum_i p(\Delta(\mathcal{D}_i))\neq0$ are in one to one correspondence with rays in $\mathfrak{D}$ with attached monomial containing $u_{l,a}$.  Such sets are naturally recovered from $\exp([W_{k,0}(A(Q'))-W_{basic}(A(Q'))]/\hbar)$.  Define $L_{P,j,a}$ to be the sum of monomials in $L_{P,j}$ that include factor $u_{l,a}$.

To find terms from $\exp([W_{k,0}(A(Q'))-W_{basic}(A(Q'))]/\hbar)$ that will contribute to $L_{P,j,a}$ upon crossing a wall that radiates from $P_l$, we should examine those not containing the factor $u_{l,w}$ for any $w$.  We consider a term $c\hbar^{-\nu} z^{\hat{n}}$ of the form:
$$c\hbar^{-\nu}z^{\hat{n}}=\hbar^{-\nu}\prod_{w=1}^\nu Mult(\mathcal{D}_w)z^{\Delta(\mathcal{D}_w)}u_{\mathcal{D}},$$
where each of the disks $\mathcal{D}_w\in RootDisk(A(Q'), r^{\mathcal{D}_w}, T_{0,tr}^{0})$.  As opposed to the case considered in \cite{kan} Lemma 5.59, we will have to consider the walls $\mathfrak{d}$ resulting from trees containing semirigid disks corresponding to unbounded rays (translated copies of $\rho_i$) emanating from $P_l$. 
Write $\hat{n}=\sum_{v=1}^\nu\Delta(\mathcal{D}_v)=\sum_{v=0}^2 n_vt_v$ and choose the primitive normal vectors $\mathbf{n}_\mathfrak{d}$ to each ray $\mathfrak{d}$ issuing from $P_l$ such that they point in the direction opposite to $\xi'_P$ when $\xi_P$ crosses $\mathfrak{d}$.  

The term $c\hbar^{-\nu}z^{\hat{n}}$ can only contribute to $L_{P,j,a}$ when $\xi_P$ crosses rays whose corresponding tree contains exactly $a+1$ semirigid disks joined at $P_l$.  The relevant rays can be enumerated as follows.  Select $\{\mathcal{D}_{i_1},\ldots \mathcal{D}_{i_s}\}\subseteq \{\mathcal{D}_1,\ldots \mathcal{D}_\nu\}$ and $M_v$ copies (here it's convenient to consider $M_v$ as an integer rather than a set) of the disk composed of the ray parallel to $\rho_v$ for $0\leq v\leq 2$ such that $s+M_0+M_1+M_2=a+1$.  Set $M=\sum_{v=0}^2 M_v t_v\in T^\Sigma$.  Let $\tilde{n}:=\sum_{v=1}^s \Delta(\mathcal{D}_{i_v})+M$ and $p(\tilde{n}):=w_{\tilde{n}}\mathbf{m}_{\tilde{n}}$, where $\mathbf{m}_{\tilde{n}}$ is primitive. These choices will produce a ray $\mathfrak{d}\in \mathfrak{D}$ with attached function 
$$f_\mathfrak{d}=1+w_{\tilde{n}}\prod_{m=1}^s Mono(\mathcal{D}_{i_m})z^{M}\frac{1}{M_0!M_1!M_2!}.$$

Let $c'h^{-(\nu-s)}z^{n'}:=\hbar^{-(\nu-s)}\prod_{\mathcal{D}\in \{\mathcal{D}_1,\ldots \mathcal{D}_\nu\}\setminus\{\mathcal{D}_{i_1},\ldots \mathcal{D}_{i_s}\}}Mono(\mathcal{D})$.  It is easy to see that the term $c'h^{-(\nu-s)}z^{n'}$ will generate a contribution of $c\hbar^{-\nu}z^{\hat{n}}$ to $L_{P,j,r}$ upon crossing $\mathfrak{d}$, and this contribution will occur exactly when $n_{j+2}+M_{j+2}\leq d =n_j+M_j<n_{j+1}+M_{j+1}$.  For simplicity of exposition, we set $j=0$ in what follows.  The quantity of the contribution is then, by definition
\begin{align*}
w_{\tilde{n}}\langle \mathbf{n}_\mathfrak{d}, \mathbf{m}_0 &\rangle D_2(d,n_0+M_0+1, n_{1}+M_1,n_{2}+M_2)\frac{\hbar^{-(\nu-s+3d-|\hat{n}|-|M|)}}{M_0!M_1!M_2!},
\end{align*}
where $D_2$ is defined in \cite{kan}, Lemma 5.43.
Noting that $|M|=a+1-s$ and recalling the isomorphism of ${\bigwedge}^2M$ with $\mathbb{Z}$, we see that the above becomes
\begin{align*}
\left(p(\tilde{n})\wedge \mathbf{m}_0 \right)D_2(d,n_0+M_0+1, n_{1}+&M_1,n_{2}+M_1)\frac{\hbar^{-(\nu-a-1+3d-|\hat{n}|)}}{M_0!M_1!M_2!}.
\end{align*}
Our goal is to now sum this contribution over all choices of $s$, $\{\mathcal{D}_{i_1},\ldots \mathcal{D}_{i_s}\}\subseteq \{\mathcal{D}_1,\ldots \mathcal{D}_\nu\}$, $M_0$, $M_1$, and $M_2$.  These should exhaust the set of relevant rays emanating from $P_l$ that $\xi_P$ crosses, and should thus calculate the total contribution.  After a little rearrangement (see \cite{thesis}), the sum becomes the following, where $t:=a+1-d+n_0$:
\begin{align}
\frac{\hbar^{-(\nu-a-1+3d-|\hat{n}|)}}{(d-n_0)!}\sum_{s=0}^{a+1}&\sum_{M_1+M_2=t-s}\frac{(-1)^{M_1+n_1+d+1}(n_1+M_1-d-1)!}{M_1!M_2!(d-n_2-M_2)!}\cdot\nonumber\\
&\left({\nu-1\choose s-1}(n_2-n_1)+{\nu\choose s}(M_2-M_1)\right).\label{mess}
\end{align}
In this sum, we are taking any terms involving factorials with negative arguments to be 0. 
\begin{sublemma}
\label{lem63}
Let $d>0$, $\nu, n_0, n_1, n_2, a\in \mathbb{Z}_{\geq 0}$ with $n_2, n_0 \leq d$. Set $t=a+1-d+n_0$, $|n|=n_0+n_1+n_2$.  Then
\begin{align*}
\frac{1}{(d-n_0)!}\sum_{s=0}^{a+1}\sum_{M_1+M_2=t-s}\frac{(-1)^{M_1+n_1+d+1}(n_1+M_1-d-1)!}{M_1!M_2!(d-n_2-M_2)!}\cdot\\
\left({\nu-1\choose s-1}(n_2-n_1)+{\nu\choose s}(M_2-M_1)\right)=\\
-\frac{1}{(d-n_0)!(d-n_1)!(d-n_2)!}{\nu+3d-|n|-1 - \left((d-n_0)+(d-n_1)\right) \choose r- \left((d-n_0)+(d-n_1)\right)}
\end{align*}
where all terms are taken to be $0$ if they involve any factorials with negative arguments.
\begin{proof}
The left hand side of the statement can be rewritten:
\begin{align*}
\frac{1}{(d-n_0)!(d-n_2)!}\sum_{M_1+M_2+s=t}&\frac{(-1)^{M_1+n_1+d+1}(n_1+M_1-d-1)!}{M_1!}{d-n_2\choose M_2}\cdot\\&\left({\nu-1\choose s-1}(n_2-n_1)+{\nu\choose s}(M_2-M_1)\right)
\end{align*}
There are two cases to distinguish: either $n_1\geq d+1$ or $n_1<d+1$.
In the former, the sum (up to a sign that won't end up mattering) in the statement should be the coefficient of $x^t$ in the expansion of:
\begin{align*}
(x+1)^{d-n_2}&x(x+1)^{\nu-1}\frac{(n_1-d-1)!}{(x+1)^{n_1-d}}(n_2-n_1)+\left(x\frac{d}{dx}\left((x+1)^{d-n_2}\right)\right)(x+1)^{\nu}\cdot\\
&\frac{(n_1-d-1)!}{(x+1)^{n_1-d}}-(x+1)^{d-n_2}(x+1)^{\nu}\left(x\frac{d}{dx}\left(\frac{(n_1-d-1)!}{(x+1)^{n_1-d}}\right)\right)
\end{align*}  
Letting $F(x)=x(x+1)^{2d-n_2-n_1+\nu-1}(n_1-d-1)!$, the above simplifies to $$(n_2-n_1)F(x)+(d-n_2)F(x)+(n_1-d)F(x)=0.$$
Note that here we didn't make use of the assumption that $d>0$.

In the latter case, we will make use of a set of functions $g_n(x)$ with $\frac{d^m}{dx^m}g_n(x)=g_{n-m}(x)$ and $g_0(x)=\frac{1}{1+x}$.  Such a set can be defined recursively as follows:
\begin{align*}
g_1(x)&=\log(1+x)\nonumber\\
g_n(x)&=\frac{(1+x)^{(n-1)}}{((n-1)!)^2}\left((n-1)!\log(1+x)-k_n\right)\nonumber\\
\end{align*}
where $k_{n+1}=(n-1)!+k_n(n)$.  Note that $g_n = \frac{n}{x+1}g_{n+1}+\frac{(x+1)^{n-1}}{(n)!}$.  By integrating, it's easy to see that the expansion of $g_n(x)$ about 0 is given by
\begin{align*}
g_n(x)=-\frac{k_nx^0}{((n-1)!)^2 0!}-\frac{k_{n-1}x^1}{((n-2)!)^21!}-\ldots-&\frac{k_{1}x^{n-1}}{((0)!)^2(n-1)!}+\\&\frac{x^n}{(n)_n}-\frac{x^{n+1}}{(n+1)_n}+\frac{x^{n+2}}{(n+2)_n}-\cdots,
\end{align*}
where $(b)_n$ is the Pochhammer symbol denoting the $n$-th falling factorial of $b$.  

Let $H_n$ denote the $n$-th harmonic number.  Note that $k_1=0$, $k_2=H_1$, and, by induction, $k_n=H_{n-1}(n-1)!$.

Let $\overbrace{f}^m$ denote the function arrived at by neglecting all terms of the expansion of $f$ with exponent less than $m$.  For example,
$$\overbrace{g_n}^n=\frac{x^n}{(n)_n}-\frac{x^{n+1}}{(n+1)_n}+\frac{x^{n+2}}{(n+2)_n}-\cdots.$$
This removes all powers of $g_n$ which are not attached to falling factorials, allowing us to write out the generating function for our sum when $d>n_1+1$.

In particular, the sum on the left hand side is the coefficient of $x^t$ in the expansion of the following generating function (about $x=0$):
\begin{align*}
P(x):=&x(x+1)^{d-n_2}(x+1)^{\nu-1}\overbrace{g_{d+1-n_1}(x)}^{d+1-n_1}(n_2-n_1)\\&+\left(x\frac{d}{dx}\left((x+1)^{d-n_2}\right)\right)(x+1)^{\nu}\overbrace{g_{d-n_1+1}(x)}^{d+1-n_1}\\
&-(x+1)^{d-n_2}(x+1)^{\nu}\overbrace{\left(x\frac{d}{dx}\left(g_{d+1-n_1}(x)\right)\right)}^{d+1-n_1}.
\end{align*}
We concentrate on the summand in the last line.
\begin{align*}
(x+1&)^{d-n_2}(x+1)^{\nu}\overbrace{\left(xg_{d-n_1}(x)\right)}^{d+1-n_1}=\\
&(x+1)^{d-n_2}(x+1)^{\nu}\overbrace{\left(x\frac{d-n_1}{x+1}g_{d+1-n_1}(x)-\frac{x(x+1)^{(d-n_1-1)}}{(d-n_1)!}\right)}^{d+1-n_1}=\\
&(x+1)^{d-n_2}(x+1)^{\nu}\overbrace{\left(x\frac{d-n_1}{x+1}g_{d+1-n_1}(x)\right)}^{d+1-n_1}-\\&\hspace{.5in}(x+1)^{d-n_2}(x+1)^{\nu}\overbrace{\left(\frac{x(x+1)^{(d-n_1-1)}}{(d-n_1)!}\right)}^{d+1-n_1}=\\
&(x+1)^{d-n_2}(x+1)^{\nu}\overbrace{\left(x\frac{d-n_1}{x+1}g_{d+1-n_1}(x)\right)}^{d+1-n_1}=\\
&(d-n_1)(x+1)^{d-n_2}(x+1)^{\nu}\overbrace{\left(x\frac{1}{x+1}g_{d+1-n_1}(x)\right)}^{d+1-n_1}\\
\end{align*}
We compare $\overbrace{\left(x\frac{1}{x+1}g_{d+1-n_1}(x)\right)}^{d+1-n_1}$ with $x\frac{1}{x+1}\overbrace{\left(g_{d+1-n_1}(x)\right)}^{d+1-n_1}$ by calculating them in terms of a product of power series expansions.  Let the coefficients of the expansion of $g_n(x)$ define $a_i, b_i\in \mathbb{Q}$ in the following way:
$$g_n(x)=\sum_{i=0}^{n-1}a_ix^i+\sum_{i=n}^\infty b_i x^i,$$
with $a_i=0$ for $i\notin \{0,\ldots, n-1\}$ and $b_i=0$ for $i\notin \{n, n+1, \ldots\}$.
\begin{align*}
\frac{1}{1+x}g_n(x)=\sum_{i=0}^\infty \left(\sum_{c=0}^{i} (a_{i-c}+b_{i-c})(-1)^{c}\right)x^i
\end{align*}
Define $\mu=a_{n-1}-a_{n-2}+\cdots+(-1)^{n-1} a_0$.  Then
\begin{align*}
\overbrace{\frac{x}{1+x}g_n(x)}^n&=\sum_{i=n-1}^\infty \left(\sum_{c=0}^{i} (a_{i-c}+b_{i-c})(-1)^{c}\right)x^{i+1}\\
&=\sum_{i=n-1}^\infty \left((-1)^{i+n-1}\mu+\sum_{c=0}^{i} (b_{i-c})(-1)^{c}\right)x^{i+1}.
\end{align*}
Alternately,
\begin{align*}
\frac{x}{1+x}\overbrace{g_n(x)}^n&=\sum_{i=0}^\infty \left(\sum_{c=0}^{i} (b_{i-c})(-1)^{c}\right)x^{i+1}\\
&=\sum_{i=n-1}^\infty \left(\sum_{c=0}^{i} (b_{i-c})(-1)^{c}\right)x^{i+1}.
\end{align*}
Thus
\begin{align*}
\overbrace{\frac{x}{1+x}g_n(x)}^n=\frac{x}{1+x}\overbrace{g_n(x)}^n+\frac{x^n}{1+x}\mu.
\end{align*}
Then the expression for $P(x)$ becomes
\begin{align*}
P(x):=&x(x+1)^{d-n_2}(x+1)^{\nu-1}\overbrace{g_{d+1-n_1}(x)}^{d+1-n_1}(n_2-n_1)\\&+x(d-n_2)(x+1)^{d-n_2-1}(x+1)^{\nu}\overbrace{g_{d+1-n_1}(x)}^{d+1-n_1} \\
&-(d-n_1)(x+1)^{d-n_2}(x+1)^{\nu}\frac{x}{x+1}\overbrace{g_{d+1-n_1}}^{d+1-n_1}\\&-(d-n_1)(x+1)^{d-n_2}(x+1)^{\nu}\frac{x^{d-n_1+1}}{x+1}\mu\\
&=-(d-n_1)(x+1)^{d-n_2}(x+1)^{\nu}\frac{x^{d-n_1+1}}{x+1}\mu
\end{align*}
To calculate $\mu$, note
\begin{align*}
a_i&=\frac{-k_{d-n_1+1-i}}{((d-n_1-i)!)^2i!}=-\frac{H_{d-n_1-i}(d-n_1-i)!}{((d-n_1-i)!)^2i!}\\&=-\frac{1}{(d-n_1)!}{d-n_1 \choose d-n_1-i}H_{d-n_1-i}.
\end{align*}
Thus, 
$$\mu=\sum_{i=0}^{d-n_1} (-1)^i a_{d-n_1-i}=\frac{-1}{(d-n_1)!}\sum_{i=0}^{d-n_1} (-1)^i H_i{d-n_1 \choose i}$$
Define $\hat{\mu}=(d-n_1)!\mu=\sum_{i=0}^{d-n_1} (-1)^i H_i{d-n_1 \choose i}$.  We claim that $\hat{\mu}=\frac{-1}{d-n_1}$.  Assuming this result, we see that this implies 
\begin{align*}
P(x)=\frac{-1}{(d-n_1)!}(x+1)^{d-n_2}(x+1)^{\nu-1}x^{d-n_1+1}.
\end{align*}  
The coefficient of $x^t$ for this final quantity is just $\frac{1}{(d-n_1)!}{\nu+d-n_2-1 \choose t-1-d+n_1}$, which gives us the desired equality.

Now for the claimed result about $\mu$.
\begin{claim}
\label{angela}
$\sum_{i=0}^{d-n_1} (-1)^i H_i{d-n_1 \choose i}=\frac{-1}{d-n_1}$
\begin{proof}

Define $B_n:=\sum_{i=0}^{n} (-1)^i H_i{n \choose i}$.  We will examine $B_n-B_{n-1}$.
\begin{align*}
B_n-B_{n-1}&=\sum_{i=0}^{n} (-1)^i H_i\left({n \choose i}-{n-1 \choose i}\right)\\
&=\sum_{i=1}^{n} (-1)^i (\frac{1}{i}+H_{i-1})\left({n-1 \choose i-1}\right)\\
&=\sum_{i=1}^{n} (-1)^i \frac{1}{i}{n-1 \choose i-1}+\sum_{i=1}^{n} (-1)^i H_{i-1}{n-1 \choose i-1}\\
&=\sum_{i=1}^{n} (-1)^i \frac{1}{n}{n \choose i}-\sum_{i=0}^{n-1} (-1)^{i} H_{i}{n-1 \choose i}\\
&=\frac{1}{n}\sum_{i=1}^{n} (-1)^i{n \choose i}-B_{n-1}\\
\end{align*}
Applying the identity $\sum_{i=0}^{n} (-1)^i{n \choose i}=0$, (unless $n=0$, which we won't be considering), we have $B_n-B_{n-1}=0-\frac{1}{n}{n\choose 0}-B_{n-1}$, which implies that $B_n=-\frac{1}{n}$.  Applying this result to our special case, we see that it proves the claim.
\end{proof}
\end{claim}

\end{proof}
\end{sublemma}
Given a non-zero contribution to  $-L_{P,j,a}$ of the term  $c\hbar^{-\nu}z^{\hat{n}}$ (with $d>0$), we can assemble a balanced tropical curve $\mathcal{C}$.   Begin by gluing the disks $\mathcal{D}_1,\ldots, \mathcal{D}_\nu$ together at their outgoing vertices at $P_l$, add on $d-n_j$ unbounded edges in the direction $\mathbf{m}_j$ for $0\leq j \leq 2$ and two additional edges $E_x$ and $E_{p}$ that will be collapsed to mark $x$ and $P_l$.  This procedure yields a frame whose valence at the new vertex $V$ is given by $Val(V):=\nu+3d-|\hat{n}| +2$.  Thus we have a tropical curve $\mathcal{C}$ with $h:\Gamma\rightarrow M_\mathbb{R}$  with $h(E_x)=h(p)=h(V)=P_l$.  Define $r:=\sum_{m=1}^\nu r^{\mathcal{D}_m}$.
The previous sublemma allows us to easily describe the contribution to $-L_{P,j,r}$ of the term $c\hbar^{-\nu}z^{\hat{n}}$ upon crossing the corresponding rays  radiating from $P_l$ as 
\begin{align*}
&{Val(V) -3 -\left((d-n_0)+(d-n_1)\right) \choose a- \left((d-n_0)+(d-n_1)\right)}\cdot\\&\hspace{1in}Mult_x^0(\mathcal{C})\left(\prod_{l=1}^\nu Mult(\mathcal{D}_l)u_{\mathcal{D}_l}\right) u_{l,a}
\hbar^{-(\nu-a-1+3d-|\hat{n}|)}
\end{align*} 
Suppose that $|r|=3d-\nu' $ for some $\nu'\geq 0$.  Let $v$ be the valence of $V$.  By construction, it is equal to $\nu+3d-|\hat{n}|+2$. On the other hand, because $\mathcal{C}$ is obtained by gluing $v-2$ semirigid disks at $V$, we have
\begin{align*}
v-2 &= \sum_{i=1}^\nu(|\Delta(\mathcal{D}_i)|-|r^{\mathcal{D}_i}|)+3d-|{\hat{n}}|\\
&=|\hat{n}|-[3d-\nu'-(a+1)]+3d-|{\hat{n}}|\\
&=\nu'+a+1
\end{align*}
Therefore, the contribution to $-L_{P,j,a}$ from $\xi_P$ crossing rays associated to this term is precisely the contribution of $\mathcal{C}$ to 
$$\langle \psi^{r(1)-1}P_{r\{1\}},\ldots, \psi^{r(\#(r))-1}P_{r\{\#(r)\}}, \psi^a P_l, \psi^{\nu'} S_2(A)\rangle^{trop}_{d,\sigma_{j,j+1}}u_{r}u_{l,a}\hbar^{-\nu'}.$$  Conversely, it is easy to see that any such curve $h$ contributing to the invariant will be accounted for by the integral by deconstructing it into its constituent semirigid disks.

Suppose $d=0$.  An examination of Expression \ref{mess} shows that any non-zero contribution must occur when $n_0=n_2=0.$  In this case, $M_2=0$, which forces $M_1=t-s=a+1-s$, so our quantity becomes
\begin{align*}
\hbar^{-(\nu-a-1-(n_{1}))}\sum_{s=0}^{a+1}\frac{(-1)^{a-s-n_1}(n_1+a-s)!}{(a+1-s)!}\left({\nu-1\choose s-1}(-n_1)+{\nu\choose s}(-a-1+s)\right).
\end{align*}
If $n_1>d=0$, then the argument applied in the first case of Lemma \ref{lem63} shows that the above quantity is equal to 0.  If $n_1=0$ then $\nu=0$, so the above simplifies to
\begin{align*}
\hbar^{-(-a-1))}\sum_{s=0}^{a+1}\frac{(-1)^{a-s-n_1}(n_1+r-s)!}{(a+1-s)!}\left({0\choose s}(-a-1+s)\right)\\
\hbar^{a+1}\frac{(-1)^{a}(a)!}{(a+1)!}\left({0\choose 0}(-a-1)\right).
\end{align*}
In this case the contribution to  $-L_{P,j,a}$ from $\xi_P$ is equal to $-(-\hbar^{a+1})u_{l,a}$.   
\end{proof}
\end{lemma}

\begin{lemma}
\label{lem64}
\begin{align*}
-L&_{i,\xi_j,\rho_{j+1}\rightarrow\sigma_{j,j+1}}^d=\hbar\sum_{P_l\in Q+\sigma_{j,j+1}}\delta_{d,0}\delta_{2,i}\left(u_{l,0}-u_{l,1}\hbar+\ldots+u_{l,k}(-\hbar)^{k}\right)+\\&\sum_{\substack{\nu\geq i-1\\ r\in \mathcal{R}_k\\|r|=3d-2+i-\nu}}\langle \psi^{r(1)-1}P_{r\{1\}},\ldots, \psi^{r(\#(r))-1}P_{r\{\#(r)\}}, \psi^\nu S_i(A) \rangle^{\rm trop}_{d,\sigma_{j,j+1}j}u_{r}h^{-(\nu+2-i)}.
\end{align*}

\begin{proof}
This follow from the previous lemmas.  Note, in particular, the first sum that results from the previous remark as $r$ is varied from $0$ to $n$.
\end{proof}
\end{lemma}
Consolidating the results of this section, we obtain the following lemma, from which Theorem \ref{yfbig} follows directly.
\begin{lemma}
\begin{align*}
L_0^d: &=\delta_{0,d}+\sum_{w\geq 0} \frac{\hbar^{-1}}{w!}\langle \frac{S_0(A)}{\hbar - \psi}, \gamma_{a,tr}^w \rangle^{trop}_{0,d}\kappa^d\\
L_1^d: &=\sum_{w\geq 0} \frac{1}{w!}\langle \frac{S_1(A)}{\hbar - \psi}, \gamma_{a,tr}^w \rangle^{trop}_{0,d}\kappa^d\\
L_2^d: &=\delta_{0,d}\hbar \sum_{j=0}^ky_{2,j}(-\hbar)^j+\sum_{w\geq 0} \frac{\hbar}{w!}\langle \frac{S_2(A)}{\hbar - \psi}, \gamma_{a,tr}^w \rangle^{trop}_{0,d}\kappa^d.
\end{align*}
\end{lemma}

\section{Formal operations} \label{form}
It should not be exceptionally difficult to use the scattering approach to directly evaluate the integrals on the potential $W_{k,\overline{m}}(A)$.  For our purposes, it is convenient to instead use the axioms of GW theory to assemble it from the result of Theorem \ref{theorem1}.  As the integral is independent of the general arrangement $A$ chosen, we will write $W_{k,\overline{m}}(A)$ as $W_{k,\overline{m}}$ in the following. 
We introduce a pair of operators on $ \mathbb{C}[T_\Sigma]\otimes_\mathbb{C} \mathfrak{R}_{k,\overline{m}}$ closely related to the \emph{fundamental class axiom} of GW theory:
$$\op:= \sum_{1\leq j, l\leq k}  u_{j,l} \frac{\partial}{\partial u_{j,l-1}}$$
and 
$$\tilop:=\exp(y_{0,0}\op)=\sum_{j=0}^\infty \frac{y_{0,0}^j}{j!}\op^j.$$
The following technical lemmas allow us to extend Theorem \ref{theorem1}.
\begin{lemma}

$$W_{k,\overline{m}}=y_{0,0}+\tilop(W_{k,0})$$ 
\begin{proof}
Let $\mathcal{D}$ be a disk in $RootDisk(A, r, T_{0,tr}^{m})$.  We will say that two disks are \emph{similar} if they differ by a permutation of the markings on the collapsed edges $E_{q_i}$.  For $1\leq j \leq k$, let $g_j$ denote the number of edges marked by elements of $\{q_1,\ldots, q_m\}$ that map to $P_j$ under $h$.  Then there are ${m \choose g_1, \ldots, g_k}$ similar disks associated to $\mathcal{D}$ which contribute a total of 
\begin{align}
\frac{y_{0,0}^m}{g_1!\ldots g_k!}Mult(\mathcal{D})u_{\mathcal{D}}z^{\Delta(\mathcal{D})}\label{eqlocal1}
\end{align}
to $W_{k,\overline{m}}$.  Define $\mathcal{D}'\in RootDisk(A, r', T_{0,tr}^{0})$ to be the result of removing the edges marked by $q_1, \ldots, q_m$ from $\mathcal{D}$ and adjusting the entries of the vector $r$ in the necessary way (removing the edges $E_{q_i}$ reduces the valencies of the vertices to which they are attached).  This disk contributes  
\begin{align}
Mult(\mathcal{D'})u_{\mathcal{D'}}z^{\Delta(\mathcal{D'})}=Mult(\mathcal{D})u_{\mathcal{D'}}z^{\Delta(\mathcal{D})}\label{eqlocal2}
\end{align}
to $W_{k,0}$.  The term $\frac{y_{0,0}^m}{m!}\op^m$ in $\tilop$ will create summands of the same multi-degree as \ref{eqlocal1} when acting on \ref{eqlocal2}, and the contribution of these terms to $\tilop(W_{k,0})$ is easily seen to equal expression  \ref{eqlocal1}.  On the other hand, it's clear how to associate a set of similar disks to any term appearing in the expansion $\tilop(W_{k,0})$ by adding marked edges $E_{q_i}$ to the associated disk in $W_{k,0}$.  Finally, the term $y_{0,0}$ in the RHS of the lemma corresponds to the semirigid disk consisting of a single $q_1$-marked edge mapping to $Q$.

\end{proof}
\end{lemma}
 \begin{lemma}
$
e^{\tilop(W_{k,0})/\hbar}=\tilop \left(e^{W_{k,0}/\hbar}\right)$
\begin{proof}
Set $r\in \mathcal{R}_k.$   
Let $\mathcal{D}_j \in  RootDisk(A, r^j, T_{0,tr}^{0})$ for $1\leq j \leq \nu$, with $r^j$ pairwise disjoint and $r$ dominating $\sum_{i=1}^j r^j:=r'$.  These disks contribute to $W_{k,0}$ and its exponential, and thus to the quantities appearing on either side of the lemma. We will compare their contribution on either side of the desired equality to terms of multi-degree $u_r y_0^{|r-r'|}\hbar^{-\nu}$.  On the LHS, this is given by 
$$\prod_{j=1}^\nu \frac{y_{0,0}^{|r-r^j|}}{|r-r^j|!}{ |r-r^j|\choose r-r^j}Mono(\mathcal{D}_j)$$  
(recall the definitions given in Section \ref{def}),  while on the RHS it is given by
$$\frac{y_{0,0}^{|r-r'|}}{|r-r'|!}{ |r-r'|\choose r-r'}\prod_{j=1}^\nu Mono(\mathcal{D}_j).$$
Because $\sum_{j=1}^\nu r-r^j=r-r'$ , the two expressions are equal.  All terms appearing on either side of the desired equality result from such choices of sets of disks, and the lemma is proven.  
\end{proof}
\end{lemma}
Together, the two previous lemmas yield:
\begin{corollary}\label{bootstrap}
$e^{W_{k,\overline{m}}/\hbar}=\tilop \left(e^{(y_{0,0}+W_{k,0})/\hbar}\right).$
\end{corollary}
If we, by abuse of notation, extend $\op$ and $\tilop$ to their obvious operators on 
$\mathbb{C}[[y_{1,0}]]\otimes_\mathbb{C} \mathfrak{R}_{k,\overline{m}}$, their actions commute  with the integration of Theorem \ref{theorem1}.  
Let $$\gamma_{b,tr}:=T_{0,tr}y_{0,0} +\gamma_{a,tr}$$ be a formal expression as before.
\begin{corollary} \label{cor1}
\begin{align*}
\sum_{i=0}^2 \alpha^i \int_{\Xi_i} e^{W_{k,\overline{m}}/\hbar} \Omega=\hbar^{-3\alpha} \sum_{j=0}^2  \left(\alpha\hbar\right)^j e^{y_{1,0}\alpha} \tilde{\Theta}_j
\end{align*}
where
\begin{align*}
\tilde{\Theta}_0: &=e^{y_{0,0}/\hbar}+\sum_{d>0,w\geq 0} \frac{\hbar^{-1}}{w!}\langle \frac{S_0(A)}{\hbar - \psi}, \gamma_{b,tr}^w \rangle^{trop}_{0,d}e^{y_{1,0}d}\\
\tilde{\Theta}_1: &=\sum_{d>0,w\geq 0} \frac{\hbar^{-1}}{w!}\langle \frac{S_1(A)}{\hbar - \psi}, \gamma_{b,tr}^w \rangle^{trop}_{0,d}e^{y_{1,0}d}\\
\tilde{\Theta}_2: &=\hbar^{-1}e^{y_{0,0}/\hbar}\sum_{j=0}^k(-\hbar)^j\sum_{l=0}^{\overline{m}} \frac{y_{0,0}^l}{l!} y_{2, l+k}+\sum_{d>0,w\geq 0} \frac{\hbar^{-1}}{w!}\langle \frac{S_2(A)}{\hbar - \psi}, \gamma_{b,tr}^w \rangle^{trop}_{0,d}e^{y_{1,0}d}
\end{align*}
in $\mathbb{C}[[y_1]]\otimes_\mathbb{C} \mathfrak{R}_{k,\overline{m}}$.
\begin{proof}
By Corollary \ref{bootstrap}
$$\sum_{i=0}^2 \alpha^i \int_{\Xi_i} e^{W_{k,\overline{m}}/\hbar} \Omega=e^{y_{0,0}/\hbar}\tilop\left(\sum_{i=0}^2 \alpha^i \int_{\Xi_i} e^{W_{k,0}/\hbar} \Omega\right).$$
Then \begin{align*}
\tilde{\Theta}_0: &=e^{y_{0,0}/\hbar}+\sum_{d>0,w\geq 0} e^{y_{0,0}/\hbar}\tilop\left(\frac{\hbar^{-1}}{w!}\langle \frac{S_0(A)}{\hbar - \psi}, \gamma_{b,tr}^w \rangle^{trop}_{0,d}e^{y_{1,0}d}\right)\\
\tilde{\Theta}_1: &=\sum_{d>0,w\geq 0} e^{y_{0,0}/\hbar}\tilop\left(\frac{\hbar^{-1}}{w!}\langle \frac{S_1(A)}{\hbar - \psi}, \gamma_{b,tr}^w \rangle^{trop}_{0,d}e^{y_{1,0}d}\right)\\
\tilde{\Theta}_2: &=e^{y_{0,0}/\hbar}\hbar^{-1}\sum_{j=0}^ky_{2,j}(-\hbar)^j+\sum_{d>0,w\geq 0} e^{y_{0,0}/\hbar}\tilop\left(\frac{\hbar^{-1}}{w!}\langle \frac{S_2(A)}{\hbar - \psi}, \gamma_{b,tr}^w \rangle^{trop}_{0,d}e^{y_{1,0}d}\right).
\end{align*}
Select $d, \nu\in \mathbb{Z}_{>0}$, $r\in \mathcal{R}_k$ with $n:=\#(r)$, and $l\in \mathbb{Z}_{\geq 0}$.  We wish to find the coefficient of $\frac{y_{0,0}^l}{l!}u_r  e^{y_{1,0}d}\hbar^{-(\nu+2)}$ in $e^{y_{0,0}/\hbar} \tilop\left(\sum_{d>0,w\geq 0} \frac{\hbar^{-1}}{w!}\langle \frac{S_0(A)}{\hbar - \psi}, \gamma_{b,tr}^w \rangle^{trop}_{0,d}e^{y_{1,0}d}\right)$.  This is readily seen to be
\begin{align*}
\sum_{i=0}^{\min(l,\nu)} \frac{l!}{i!}\sum_{\substack{r'\prec r\\|r-r'|=l-i}}\frac{1}{|r-r'|!} { |r-r'| \choose r-r' } \langle  \psi^{r'(1)-1}P_{r\{1\}},\ldots,\psi^{r'(n)-1} P_{r\{n\}},\psi^{\nu-i} S_0(A)\rangle^{trop}_{0,d} \\
=\sum_{i+w_1+\ldots+w_n=l}{ l \choose i, w_1, \ldots, w_n } \langle  \psi^{r(1)-w_1-1}P_{r\{1\}},\ldots,\psi^{r(n)-w_n-1} P_{r\{n\}},\psi^{\nu-i} S_0(A)\rangle^{trop}_{0,d},
\end{align*}
where the above invariants are interpreted as zero if they contain any negative powers of $\psi$.  By iterating the tropical fundamental class axiom (Lemma \ref{tropfun}), it's easy to see that
\begin{align*}
 \langle  \psi^{r(1)-1}P_{r\{1\}}\ldots&,\psi^{r(n)-1} P_{r\{n\}}, T_{0,tr}^l,\psi^{\nu}S_0(A)\rangle^{trop}_{0,d}
 \end{align*}
is equal to the above expression.
Of course, the same result holds when replacing $S_0(A)$ with $S_1(A)$ or $S_2(A)$.  
\end{proof}
\end{corollary}

Next, we normalize the integral from the above corollary to satisfy the conditions of Section 1 of \cite{MSP2}, allowing us to apply mirror symmetry.

\begin{lemma}\label{lemmaop}
Let $\Xi\in H_2(\kappa^{-1}(u), {\rm Re}(W_{basic}/\hbar)\ll 0, \mathbb{C})$.  Then 
\begin{align*}
\int_\Xi e^{W_{k,\overline{m}}/\hbar} (\op(W_{k,\overline{m}}))\Omega=\hbar \op\left(\int_\Xi e^{W_{k,\overline{m}}/\hbar}\Omega\right)
\end{align*}
in $\mathbb{C}[T_\Sigma]\otimes_\mathbb{C}\mathfrak{R}_{k,\overline{m}}.$
\begin{proof}
\begin{align*}
\int_\Xi e^{W_{k,\overline{m}}/\hbar} (\op(W_{k,\overline{m}}))\Omega&=\hbar \int_\Xi \op(e^{W_{k,\overline{m}}/\hbar})\Omega\\
&=\hbar \op\left(\int_\Xi e^{W_{k,\overline{m}}/\hbar}\Omega\right).
\end{align*}
\end{proof}
\end{lemma}
Combining the result of Lemma \ref{lemmaop} and the tropical fundamental class axiom (Lemma \ref{tropfun}), we immediately achieve the following result.
\begin{corollary} \label{bigcor} Let $f:=1+\op(W_{k,\overline{m}})$.  Then
\begin{align*}
\sum_{i=0}^2 \alpha^i \int_{\Xi_i} e^{W_{k,\overline{m}}/\hbar} f\Omega=\hbar^{-3\alpha} \sum_{j=0}^2  \left(\alpha\hbar\right)^j e^{y_{1,0}\alpha} L_j,
\end{align*}
where
\begin{align*}
L_0: &=e^{y_{0,0}/\hbar}+\sum_{d>0,w\geq 0} \frac{1}{w!}\langle \frac{Q}{\hbar - \psi}, T_{0,tr}, \gamma_{b,tr}^w \rangle^{trop}_{0,d}e^{y_{1,0}d}\\
L_1: &=\sum_{d>0,w\geq 0} \frac{1}{w!}\langle \frac{L}{\hbar - \psi}, T_{0,tr}, \gamma_{b,tr}^w \rangle^{trop}_{0,d}e^{y_{1,0}d}\\
L_2: &=\hbar^{-1}e^{y_{0,0}/\hbar}\sum_{l=0}^{\overline{m}} \frac{y_{0,0}^l}{l!} y_{2, l} +\sum_{d>0,w\geq 0} \frac{1}{w!}\langle \frac{M_\mathbb{R}}{\hbar - \psi}, T_{0,tr}, \gamma_{b,tr}^w \rangle^{trop}_{0,d}e^{y_{1,0}d}.
\end{align*}
If we define $\phi_i$ by rewriting
$$\hbar^{-3\alpha} \sum_{j=0}^2  \left(\alpha\hbar\right)^j e^{y_{1,0}\alpha} L_j=\hbar^{-3\alpha} \sum_{j=0}^2  \left(\alpha\hbar\right)^j \phi_j,$$
we see that
\begin{align*}
\phi_0: &=L_0\\
\phi_1: &=y_{1,0}\hbar^{-1}L_0+L_1\\
\phi_2: &=\frac{y_{1,0}^2\hbar^{-2}}{2}L_0+y_{1,0}\hbar^{-1}L_1+L_2.
\end{align*}
If we write $\phi_i:=\sum_{j=0}^\infty \hbar^{-j} \phi_{i,j}$ with $\phi_{i,j}\in \mathbb{C}[[y_{1,0}]]\otimes_\mathbb{C} \mathfrak{R}_{k,\overline{m}},$
\begin{align*}
\phi_{i,0}&=\delta_{i,0}\\
\phi_{0,1}&=y_{0,0}+\tilde{K}_2\\
\phi_{1,1}&=y_{1,0}+\tilde{K}_1\\
\phi_{2,1}&=\sum_{l=0}^{\overline{m}} \frac{y_{0,0}^l}{l!} y_{2,l}+\tilde{K}_0,
\end{align*}
where 
$$\tilde{K}_i:=\sum_{d>0,w\geq 0} \frac{1}{m!}\langle S_{2-i}(A), T_{0,tr}, \gamma_{b,tr}^w \rangle^{trop}_{0,d}e^{y_{1,0}d}.$$
\end{corollary}

By Theorem \ref{theorem1}, the expressions above are independent of the choice of arrangement $A$.   Thus, we can simply write 
$$\langle \psi^{a_1}T_{2,tr}\ldots,\psi^{a_n} T_{2,tr}, T_{0,tr}^m,\psi^\nu T_{2-i}\rangle_{0,d}^{trop}$$
 in place of 
 $$\langle \psi^{a_1}P_{r\{1\}}\ldots,\psi^{a_n} P_{r\{n\}}, T_{0,tr}^m,\psi^\nu S_i(A)\rangle_{0,d}^{trop}.$$
With this observation, we write
$$\gamma_{b,tr}:=T_{0,tr}y_{0,0}+T_{2,tr} y_{2,0}+\psi T_{2,tr} y_{2,1}+\ldots +\psi^{k-1} T_{2,tr} y_{2,k-1}.$$
Let 
$$\gamma_{b,cl}:=T_{0}y_{0,0}+T_2 y_{2,0}+\psi T_2 y_{2,1}+\ldots +\psi^{k-1} T_2 y_{2,k-1}$$
and
$$\gamma_{c,cl}:=T_{0}y_{0,0}+T_1y_{1,0}+T_2 y_{2,0}+\psi T_2 y_{2,1}+\ldots +\psi^{k-1} T_2 y_{2,k-1}$$
where $T_i$ is a positive generator of  $H^{2i}(\mathbb{P}^2, \mathbb{Z})$.
Given the results of Markwig and Rau \cite{M+R} showing the classical relevance of these tropical invariants, we have 
\begin{align*}
\tilde{K}_i&=\sum_{d>0,w\geq 0} \frac{1}{w!}\langle T_{2-i, tr}, T_{0,tr}, \gamma_{b,tr}^w \rangle^{trop}_{0,d}e^{y_{1,0}d}\\
&=\sum_{d>0,w\geq 0} \frac{1}{w!}\langle T_{2-i}, T_0, \gamma_{b,cl}^w \rangle^{cl}_{0,d}e^{y_{1,0}d}\\
&=\sum_{d>0,w\geq 0} \frac{1}{w!}\langle T_{2-i}, T_0, \gamma_{c,cl}^w \rangle^{cl}_{0,d},
\end{align*}
where we use the divisor axiom and the convention that GW invariants of incompatible dimension (see \cite{c+k} for a discussion of dimensions) are equal to zero.  Let 
$$K_i:= \sum_{d,w\geq 0} \frac{1}{w!}\langle T_{i}, T_0, \gamma_{c,cl}^w \rangle^{cl}_{0,d}$$
(we now include degree 0 invariants).

Using the fundamental class and point mapping axioms of classical GW theory to analyze the degree 0 pieces of $K_i$, one can see that $\phi_{i,1}=K_{2-i}$.  For example, 
\begin{align*}
 K_{2}&=\sum_{w\geq 0} \frac{1}{w!}\langle T_2, T_0, \gamma_{c,cl}^w \rangle^{cl}_{0,0}\\&= \langle T_2, T_0, \gamma_{c,cl} \rangle^{cl}_{0,0}\\
&= \langle T_2, T_0, y_{0,0}T_0 \rangle^{cl}_{0,0}\\
 &=  y_{0,0}
 \end{align*}
 
 We take the inverse limit 
$$\mathfrak{R}_{k}:= R_k[[y_{0,0}]]=\lim_{\leftarrow \overline{m}} \mathfrak{R}_{k,\overline{m}}$$
and extend the above results in the obvious way.
\section{Mirror Symmetry}

\begin{definition}\label{jfunction}
We consider Givental's $J$ function as an element $$J_{\mathbb{P}^2} \in  \mathbb{C}[[\tilde{y_0}, \tilde{y_{1}}, \tilde{y_2}, \hbar^{-1}]]\otimes H^*(\mathbb{P}^2, \mathbb{Z}),$$ defined as in \cite{iri}, up to some minor rearrangement, as
\begin{align*}
J_{\mathbb{P}^2}=&e^{(T_0 \tilde{y_0}+T_1 \tilde{y_1})/\hbar}\cup\Big(T_0+\tilde{y_2}T_2+
\\&\sum_{i=0}^2\big(\sum_{d\geq 1,\nu\geq 0} \langle T_2^{3d+i-2-\nu}, \psi^\nu T_{2-i}\rangle^{cl}_{0,d}\hbar^{-(\nu+2)}e^{d\tilde{y_{1}}}\frac{\tilde{y_2}^{3d+i-2-\nu}}{(3d+i-2-\nu)!}\big)T_i\Big).
\end{align*}
Define $J_i$ to be the $T_i$ component of $J$.
\end{definition}
\begin{lemma}\label{jfun}
Let $\gamma:=T_0 \tilde{y_0}+T_1 \tilde{y_1}+T_2 \tilde{y_2}$
\begin{align*}
J_{\mathbb{P}^2}=T_0+\sum_{w,d=0}^\infty \sum_{i=0}^2 \frac{1}{w!} \langle \frac{T_{2-i}}{\hbar-\psi}, T_0, \gamma^w\rangle^{cl}_{0,d} T_i
\end{align*}
\begin{proof}
Let us examine $J_1$.  
\begin{align}
J_1=&\frac{\tilde{y}_1 e^{\tilde{y}_0/\hbar}}{\hbar}\left(1+\sum_{d\geq 1}\sum_{j,\nu\geq 0} \langle T_2^{j}, \psi^\nu T_{2}\rangle^{cl}_{0,d}\hbar^{-(\nu+2)}e^{d\tilde{y_{1}}}\frac{\tilde{y_2}^{j}}{j!}\right)\nonumber\\+
&e^{\tilde{y}_0/\hbar}\left(\sum_{d\geq 1}\sum_{j,\nu\geq 0} \langle T_2^{j}, \psi^\nu T_{1}\rangle^{cl}_{0,d}\hbar^{-(\nu+2)}e^{d\tilde{y_{1}}}\frac{\tilde{y_2}^{j}}{j!}\right)\nonumber\\
=&\frac{\tilde{y}_1 }{\hbar}\left(e^{\tilde{y}_0/\hbar}+\sum_{d\geq 1}\sum_{j,w,\nu\geq 0} \langle T_2^j, T_0^w,\psi^{\nu+w} T_{2}\rangle^{cl}_{0,d}\hbar^{-(w+\nu+2)}e^{d\tilde{y_{1}}}\frac{\tilde{y}_0^w\tilde{y_2}^{j}}{w!j!}\right)\label{line1}\\+
&\sum_{d\geq 1}\sum_{j,w, \nu\geq 0} \langle T_2^j,  T_0^w,\psi^{\nu +w}T_{1}\rangle^{cl}_{0,d}\hbar^{-(\nu+w+2)}e^{d\tilde{y_{1}}}\frac{\tilde{y}_0^w\tilde{y_2}^j}{w!j!}\nonumber\\
=&\frac{\tilde{y}_1 }{\hbar}e^{\tilde{y}_0/\hbar}+\sum_{d,\nu\geq 1}\sum_{j,w\geq 0} \langle T_2^{j}, T_0^w,\psi^{\nu+w-1} T_{2}\rangle^{cl}_{0,d}\hbar^{-(w+\nu+2)}\tilde{y}_1 e^{d\tilde{y_{1}}}\frac{\tilde{y}_0^w\tilde{y}_2^j}{w!j!}\nonumber\\+
&\sum_{d\geq 1}\sum_{j,w,\nu\geq 0} \langle T_2^j,  T_0^w,\psi^{\nu+w}T_{1}\rangle^{cl}_{0,d}\hbar^{-(w+\nu+2)}e^{d\tilde{y_{1}}}\frac{\tilde{y}_0^w\tilde{y_2}^j}{w!j!}\nonumber
\end{align}
\begin{align}
=&\sum_{t,\nu \geq 0} \frac{1}{t!}\langle \psi^\nu T_1, T_0, \gamma^t \rangle^{cl}_{0,0}+\label{line3}\\&\sum_{d,l,\nu\geq 1}\sum_{j,w\geq 0} \langle T_2^j, T_0^w,\psi^{\nu+w-1} T_{2}\rangle^{cl}_{0,d}\hbar^{-(w+\nu+2)}l\frac{\tilde{y}_0^wd^{l-1}\tilde{y}_1^{l}\tilde{y}_2^j}{w!l!j!}\nonumber\\+
&\sum_{d,l\geq 1}\sum_{j, w, \nu\geq 0} \langle T_2^j,  T_0^w,\psi^{\nu+w}T_{1}\rangle^{cl}_{0,d}\hbar^{-(w+\nu+2)}\frac{\tilde{y}_0^wd^l\tilde{y}_1^l\tilde{y_2}^j}{w!l!j!}\nonumber\\
=&\sum_{t,\nu \geq 0} \frac{1}{t!}\langle \psi^\nu T_1, T_0, \gamma^t \rangle^{cl}_{0,0}+\label{line4}\\&\sum_{d, l\geq 1}\sum_{l,j,w,\nu\geq 0} \langle T_2^j,  T_0^w, T_1^l,\psi^{\nu+w}T_{1}\rangle^{cl}_{0,d}\hbar^{-(w+\nu+2)}\frac{\tilde{y}_0^w\tilde{y}_1^l\tilde{y_2}^{j}}{w!l!j!}\nonumber\\
=&\sum_{d,l,j,w,\nu\geq 0} \langle T_2^j,  T_0^w, T_1^l,\psi^{\nu+w}T_{1}\rangle^{cl}_{0,d}\hbar^{-(w+\nu+2)}\frac{\tilde{y}_0^w\tilde{y}_1^l\tilde{y_2}^j}{w!l!j!}\nonumber\\
=&\sum_{d, w,\geq 0} \frac{1}{w!}\langle \gamma^w, T_0,\frac{T_{1}}{\hbar-\psi}\rangle^{cl}_{0,d}\label{line5}
\end{align}
Equality \ref{line1} is due to the fundamental class axiom, \ref{line3} is due to the point mapping axiom, \ref{line4} is due to the divisor axiom, and \ref{line5} is due again to the fundamental class axiom.  The other pieces of the lemma follow from similar analysis.
\end{proof}
\end{lemma}

Define a map
$$\Phi: \mathbb{C}[[\tilde{y_0}, \tilde{y_{1}}, \tilde{y_2},\hbar^{-1}]]\otimes H^*(\mathbb{P}^2, \mathbb{Z})\rightarrow \mathfrak{R}_k [[y_{1,0}, \hbar^{-1}]]\otimes H^*(\mathbb{P}^2, \mathbb{Z})$$
by $\tilde{y}_i\mapsto K_{2-i}$ for $0\leq i \leq 2$.

Let $$\mathbb{T}_{trop}:=\sum_{i=0}^2 \phi_i T_i$$
and 
$$\mathbb{J}:=\Phi(J).$$  

\begin{theorem} 
\label{theorem2} 
Let $\mathfrak{M}_{\Sigma, k}$ be the formal spectrum of the completion of $\mathbb{C}[K_\Sigma]\otimes_\mathbb{C}\mathfrak{R}_k$ at the maximal ideal $(y_{0,0}, \kappa-1, \{u_{i,j}\}_{i,j})$.  The completion is isomorphic to $\mathbb{C}[[y_{1,0}]]\otimes_\mathbb{C} \mathfrak{R}_{k}$ with $y_{1,0}:=\log\kappa$, the latter expanded in a power series at $\kappa=1$.  Let
$$\check{\mathfrak{X}}_{\Sigma, k}=\check{\mathcal{X}}_{\Sigma, k}\times _{\mathcal{M}_{\Sigma, k}} \mathfrak{M}_{\Sigma, k}.$$
The function $W_{k,\overline{m}}$ is regular (for all $m$) on $\mathfrak{X}_{\Sigma, k}$ and restricts to $W_{basic}=x_0+x_1+x_2$ on the closed fiber of $\check{\mathfrak{X}}_{\Sigma, k}\rightarrow  \mathfrak{M}_{\Sigma, k}$ and hence gives a deformation of this function over $ \mathfrak{M}_{\Sigma, k}$.  Thus we have a morphism $\omega$ from  $\mathfrak{M}_{\Sigma, k}$ to the universal unfolding space $\mathcal{M}:=\Spec \mathbb{C}[[\tilde{y_0}, \tilde{y_{1}}, \tilde{y_2}]]$.  This map is given by:
\begin{align*}
\tilde{y_0}\mapsto K_2\\
\tilde{y_{1}}\mapsto K_1\\
\tilde{y_2}\mapsto K_0.
\end{align*}
The morphism $\omega$ induces the map $\Phi$ defined above, and 
$$\mathbb{T}_{trop}=\mathbb{J}.$$
\begin{proof}
Follows from Corollary \ref{cor1} as our data satisfies the conditions of Section 1 of \cite{MSP2}.  See Corollary 3.9 of \cite{MSP2}.
\end{proof}
\end{theorem}
For convenience, we again take an inverse limit
$$\mathfrak{R}:=\lim_{\leftarrow k} \mathfrak{R}_{k}$$
and extend the above definitions and results.  Note the natural inclusion of 
$$\tilde{\mathfrak{R}}:=\mathbb{C}[[\hbar^{-1}, y_{0,0}, y_{1,0}, y_{2,0}, y_{2,1},\ldots]]$$
 into $\mathfrak{R}[[y_{1,0}, \hbar^{-1}]]$ given by 
$$y_{2,i}\mapsto \sum_{j} u_{j,i}.$$
Because the period integrals are symmetric with respect to point labelings, the limits of $\mathbb{T}_{trop}$ and $\mathbb{J}$ in $\mathfrak{R}[[y_{1,0}, \hbar^{-1}]]$ are in the image of this inclusion.  In the following, we restrict to this setting.

Consider the classical version $\mathbb{T}$ of $\mathbb{T}_{trop}$, where each tropical invariant is replaced by its corresponding classical GW invariant.  We will show that $\mathbb{T}=\mathbb{J}$, thus implying $\mathbb{T}_{trop}=\mathbb{T}.$  Consider a component of $\mathbb{T}=\sum_{i=0}^2 \mathbb{T}_iT_i$. 
\begin{align*}
\mathbb{T}_0&=e^{y_{0,0}/\hbar}+\sum_{d>0,w\geq 0} \frac{1}{w!}\langle \frac{T_{2}}{\hbar - \psi}, T_0, \gamma_{b,cl}^w \rangle^{cl}_{0,d}e^{y_{1,0}d}\\
&=T_0+\sum_{v\geq 1} \langle \psi^{v-1} T_2, T_0, T_0^{v}\rangle^{cl}_{0,0}\frac{y_{0,0}^v}{v!\hbar^v}+\sum_{d>0,w\geq 0} \frac{1}{w!}\langle \frac{T_{2}}{\hbar - \psi}, T_0, \gamma_{b,cl}^w \rangle^{cl}_{0,d}e^{y_{1,0}d}\\
&=T_0+\sum_{d,w\geq 0} \frac{1}{w!}\langle \frac{T_2}{\hbar - \psi}, T_0, \gamma_{b,cl}^w \rangle^{cl}_{0,d}e^{y_{1,0}d}\\
&=T_0+\sum_{d,w\geq 0} \frac{1}{w!}\langle \frac{T_2}{\hbar - \psi}, T_0,\gamma_{c,cl}^w \rangle^{cl}_{0,d}\\
\end{align*}
The second equality is by the fundamental class and point mapping axioms, while the last follows from the the divisor axiom.  A parallel analysis (similar to that found in Lemma \ref{jfun}) can be applied to the other components, yielding the following identity:
\begin{align*}
\mathbb{T}=T_0+\sum_{i=0}^2 \sum_{d,w\geq 0} \frac{1}{w!} \langle \frac{T_{2-i}}{\hbar-\psi}, T_0, \gamma_{c,cl}^w \rangle^{cl}_{0,d}T_i.
\end{align*}
We now turn our attention to $\mathbb{J}$.  Defining $\gamma_{J}=T_0K_2+T_1K_1+T_2K_0$, we can write 
$$\mathbb{J}=T_0+\sum_{i=0}^2 \sum_{d,w\geq 0}\frac{1}{w!} \langle \frac{T_{2-i}}{\hbar-\psi}, T_0, \gamma_{J}^w\rangle^{cl}_{0,d} T_i,$$ a generating function whose coefficients can be written entirely in terms of the classical GW invariants of $\mathbb{P}^2$.

The following operators will be, for our purposes, closely related to the dilaton axiom. 
$$\diff:= y_{0,0} \frac{\partial}{\partial{y_{0,0}}}+y_{1,0}\frac{\partial}{\partial{y_{1,0}}}+\sum_{i\geq 0}y_{2,i}\frac{\partial}{\partial{y_{2,i}}}$$
$$\diff_>:=\sum_{i>0} y_{2,i} \frac{\partial}{\partial{y_{2,i}}}.$$
\begin{lemma}
\label{lemmaT}
\begin{align*}
\sum_{j=0}^2 \left (\frac{\partial}{\partial y_{j,0}} \mathbb{T} \right) K_{2-j}&= \diff(\mathbb{T}).
\end{align*}
\begin{proof}
\begin{align*}
\sum_{j=0}^2 &\left(\frac{\partial}{\partial y_{j,0}} \mathbb{T} \right)K_{2-j}\\
&=\sum_{j,i=0}^2 \left(\sum_{d,w\geq 0} \frac{1}{w!} \langle \frac{T_{2-i}}{\hbar-\psi}, T_0, \gamma_{c,cl}^w, T_j \rangle^{cl}_{0,d}\sum_{d',w'\geq 0} \frac{1}{w'!}\langle T_{2-j}, T_0, \gamma_{c,cl}^{w'} \rangle^{cl}_{0,d'}\right)T_i\\
&=\sum_{i=0}^2 \sum_{d,w\geq 0} \frac{1}{w!} \langle \frac{T_{2-i}}{\hbar-\psi}, T_0, \gamma_{c,cl}^w,  \psi T_0 \rangle^{cl}_{0,d}T_i\\
&=\sum_{i=0}^2 \diff\left(\sum_{d,w\geq 0} \frac{1}{w!} \langle \frac{T_{2-i}}{\hbar-\psi}, T_0, \gamma_{c,cl}^w \rangle^{cl}_{0,d}T_i\right)\\
&=\diff(\mathbb{T}).
\end{align*}
The second equality is due to the topological recursion relationship (see \cite{kan}, Proposition 2.12) , while the third is due to the dilaton axiom.
\end{proof}
\end{lemma}

\begin{lemma}
\begin{align*}
\sum_{j=0}^2 \left (\frac{\partial}{\partial y_{j,0}} \mathbb{J} \right) K_{2-j}&= \diff(\mathbb{J}).
\end{align*}
\begin{proof}
\begin{align*}
\sum_{j=0}^2& \left (\frac{\partial}{\partial y_{j,0}} \mathbb{J} \right)K_{2-j}=\sum_{j,i=0}^2 \frac{\partial}{\partial y_{j,0}}\left(\sum_{w,d=0}^\infty \frac{1}{w!} \langle \frac{T_{2-i}}{\hbar-\psi}, T_0, \gamma_{J}^w\rangle^{cl}_{0,d}\right) K_{2-j}T_i=\\
&\sum_{i=0}^2 \sum_{w\geq 1,d \geq 0} \frac{1}{(w-1)!} \langle \frac{T_{2-i}}{\hbar-\psi}, T_0, \gamma_{J}^{w-1}\rangle^{cl}_{0,d}\sum_{j=0}^2\left(\frac{\partial}{\partial y_{j,0}}\left(K_0+K_1+K_2\right)\right) K_{2-j}T_i.
\end{align*}
For $0\leq l \leq 2$,
\begin{align*}
\sum_{j=0}^2\frac{\partial}{\partial y_{j,0}}\left(K_l\right) K_{2-j}&=\sum_{d,w\geq 0} \frac{1}{w!} \langle T_l, T_0, \gamma_{c,cl}^w, T_j \rangle^{cl}_{0,d}\sum_{d',w'\geq 0} \frac{1}{w'!}\langle T_{2-j}, T_0, \gamma_{c,cl}^{w'} \rangle^{cl}_{0,d'}\\
&=\sum_{d,w\geq 0} \frac{1}{w!} \langle T_l, T_0, \gamma_{c,cl}^w,  \psi T_0 \rangle^{cl}_{0,d}\\
&=\diff(K_l),
\end{align*}
where the above equalities follow from the reasoning used in the previous lemma.  So we have 
\begin{align*}
\sum_{j=0}^2 &\left (\frac{\partial}{\partial y_{j,0}} \mathbb{J} \right)K_{2-j}\\
&=\sum_{i=0}^2 \sum_{w\geq 1,d \geq 0}^\infty \frac{1}{(w-1)!} \langle \frac{T_{2-i}}{\hbar-\psi}, T_0, \gamma_{J}^{w-1}\rangle^{cl}_{0,d}\diff\left(K_0+K_1+K_2\right)T_i\\
&=\diff\left(\sum_{i=0}^2 \sum_{w,d\geq 0}^\infty \frac{1}{(w)!} \langle \frac{T_{2-i}}{\hbar-\psi}, T_0, \gamma_{J}^{w}\rangle^{cl}_{0,d}T_i\right).
\end{align*}
\end{proof}
\end{lemma}
We use induction to show $\mathbb{T}=\mathbb{J}$.  Define a $\mathbb{Z}\left[\frac{1}{3}\right]$ grading on the monomials of $\tilde{\mathfrak{R}}\otimes H^*(\mathbb{P}^2, \mathbb{Z})$ by 
$$\gr(y_{0,0}^j y_{1,0}^l \hbar^{-\nu}\prod_{m} y_{2,m}^{a_m}T_i):= \frac{1}{3}\left(\nu -j - i + \sum_{m} a_{m}(m+1)\right)+\sum_{m>0}a_{m}.$$
Note that $\diff$ and $\diff_>$ preserve the $\gr$-grading  of monomials not sent to $0$.  

When applied to a summand of $\mathbb{T}$, 
\begin{align*}
\gr\left( \langle T_0, T_0^{r_{0,0}}, T_1^{r_{1,0}}, T_2^{r_{2,0}}, \ldots, (\psi^{k-1} T_2)^{r_{2,k-1}}, \psi^{\nu} T_{2-i}\rangle^{cl}_{0,d}\hbar^{-(\nu+1)}T_i \prod_{a,b} \frac{y_{a,b}^{r_{a,b}}}{r_{a,b}!}\right)=\\
d+\text{ the number of insertions with positive exponent on }\psi,
\\\text{ excluding the term whose power of } \psi\text{ is recorded by the exponent of }\hbar.\\
\end{align*}
The integrality results from dimension restrictions of nonzero invariants.
Defining $\hat{K}_i:=K_i T_0$, we examine the value of $\gr$ on a typical term.
\begin{align*}
\gr\left(\langle  T_0, T_0^{r_{0,0}}, T_1^{r_{1,0}}, T_2^{r_{2,0}}, (\psi^1 T_2)^{r_{2,1}}, \ldots, (\psi^{k-1} T_2)^{r_{2,k-1}},T_i\rangle^{cl}_{0,d} T_0\prod_{a,b} \frac{y_{a,b}^{r_{a,b}}}{r_{a,b}!}\right)=\\
d+\frac{1-i}{3}+\text{ the number of insertions with positive exponent on }\psi.
\end{align*}
When applied to $\mathbb{J}$, $\gr$ admits a similar description.  The coefficient of the monomial of a particular degree is a sum of products of GW invariants.  When expressed in this form, the grading can be recovered from any summand of the coefficient as the sum of the degrees of the invariants in the product with the count of the total number of insertions with non-trivial $\psi$-classes, again excluding the term whose exponent of $\psi$ is recorded by the power of $\hbar$.  

For $j\in \mathbb{Z}\left[\frac{1}{3}\right]$, define $\mathbb{J}_{[j]}$ and $\mathbb{T}_{[j]}$ to be the $\gr$-degree $j$ monomials of $\mathbb{J}$ and $\mathbb{K}$, respectively.
Note that $\gr$ is integral and non-negative for all non-zero terms in $\mathbb{J}$ and $\mathbb{T}$.  The base case $\mathbb{J}_0=\mathbb{T}_0$ follows from the point mapping axiom.  Let $n>0\in \mathbb{Z}$, and assume $\mathbb{J}_{[j]}=\mathbb{T}_{[j]}$ for all $j<n$.

We analyze the degree $n$ part of $\diff(\mathbb{T})$ using Lemma \ref{lemmaT}.
\begin{align}
\diff\left(\mathbb{T}\right)_{[n]}&=\left(\sum_{j=0}^2 \left(\frac{\partial}{\partial {y_{j,0}}}\mathbb{T}\right)\hat{K}_{2-j}\right)_{[n]}\nonumber\\
&=\sum_{w=0}^n\sum_{j=0}^2 \left(\frac{\partial}{\partial {y_{j,0}}}\mathbb{T}\right)_{[w+\frac{1-j}{3}]}\left(\hat{K}_{2-j}\right)_{[n-w-\frac{1-j}{3}]}\label{crux}\\
&=\sum_{w=0}^n\sum_{j=0}^2 \left(\frac{\partial}{\partial {y_{j,0}}}\mathbb{T}_{[w]}\right)\left(\hat{K}_{2-j}\right)_{[n-w-\frac{1-j}{3}]}\nonumber\\
&=\sum_{j=0}^2\left( \sum_{w=0}^{n-1} \left(\frac{\partial}{\partial {y_{j,0}}}\mathbb{T}_{[w]}\right)\left(\hat{K}_{2-j}\right)_{[n-w-\frac{1-j}{3}]}\right)\nonumber\\
&\hspace{1in}+\left(\frac{\partial}{\partial {y_{j,0}}}\mathbb{T}_{[n]}\right)\left(\hat{K}_{2-j}\right)_{[-\frac{1-j}{3}]}\nonumber\\
&=\sum_{j=0}^2\left( \sum_{w=0}^{n-1} \left(\frac{\partial}{\partial {y_{j,0}}}\mathbb{J}_{[w]}\right)\left(\hat{K}_{2-j}\right)_{[n-w-\frac{1-j}{3}]}\right)+\left(\frac{\partial}{\partial {y_{j,0}}}\mathbb{T}_{[n]}\right)y_{j,0}\nonumber
\end{align}
The indices in equality \ref{crux} are due to the integrality of $\gr$ on monomials of $\mathbb{T}$ and the action of $\frac{\partial}{\partial {y_{j,0}}}$.  Noting that the second summand of the last line is precisely the difference between  $\diff_> (\mathbb{T})_{[n]}$ and $\diff (\mathbb{T})_{[n]}$, 
\begin{align*}
\diff_> (\mathbb{T})_{[n]}=&\sum_{j=0}^2\left( \sum_{w=0}^{n-1} \left(\frac{\partial}{\partial {y_{j,0}}}\mathbb{J}_{[w]}\right)\left(\hat{K}_{2-j}\right)_{[n-w-\frac{1-j}{3}]}\right)\\
=&\sum_{j=0}^2\Big( \sum_{w=0}^{n} \left(\frac{\partial}{\partial {y_{j,0}}}\mathbb{J}_{[w]}\right)\left(\hat{K}_{2-j}\right)_{[n-w-\frac{1-j}{3}]}+\\
&-\left(\frac{\partial}{\partial {y_{j,0}}}\mathbb{J}_{[n]}\right)\left(\hat{K}_{2-j}\right)_{[-\frac{1-j}{3}]}\Big)\\
=&\sum_{j=0}^2\left( \left(\frac{\partial}{\partial {y_{j,0}}}\mathbb{J}\right)\left(\hat{K}_{2-j}\right)\right)_{[n]}-y_{0,j}\frac{\partial}{\partial {y_{j,0}}}\mathbb{J}_{[n]}\\
=&\diff_>({\mathbb{J}})_{[n]}.
\end{align*}

Therefore, $\diff_>(\mathbb{T})=\diff_>(\mathbb{J})$.  Of course, the above equality implies $\mathbb{T}=\mathbb{J}$ in all degrees except for those which are in the kernel of $\diff_>$, i.e. degree $0$ in $y_{2,j}$ for all $j>0$.  However, these exceptional degrees are exactly those for which Gross has shown $\mathbb{T}$ is equal to $\mathbb{T}_{trop}=\mathbb{J}.$  Therefore, we have proven the following:
\begin{theorem}\label{theorem3}
$\mathbb{T}=\mathbb{J}$.
\end{theorem}

\begin{corollary}
$\mathbb{T}_{trop}=\mathbb{T}.$  That is, the tropical descendent GW invariants of Definition \ref{tropmult} are equal to their classical counterparts.
\end{corollary}

\bibliography{bib}
\end{document}